\documentclass[12pt]{amsart}     

\usepackage{tkz-fct}

\usepackage[margin=1.25in]{geometry}

\usepackage[pagebackref=true,
]{hyperref} 
\hypersetup{colorlinks=true,linkcolor=blue,filecolor=magenta,citecolor=green,urlcolor=blue,final}

\usepackage{amsthm,amssymb,amsmath,amscd}
\usepackage{amsfonts,color}
\usepackage{latexsym}
\usepackage[all]{xy}
\usepackage{ulem}

\usepackage{mathptmx}

\newtheorem{thm}{Theorem}[section]
\newtheorem{lem}[thm]{Lemma}
\newtheorem{cor}[thm]{Corollary}
\newtheorem{prop}[thm]{Proposition}
\theoremstyle{definition}
\newtheorem{df}[thm]{Definition}
\theoremstyle{remark}
\newtheorem{rem}[thm]{Remark}

\newtheorem{example}[thm]{Example}

\numberwithin{equation}{section}

\newtheorem*{ack}{Acknowledgements}

\newcommand{\bC}{{\mathbb C}}
\newcommand{\bD}{{\mathbb D}}

\newcommand{\bL}{{\mathbb L}}

\newcommand{\bR}{{\mathbb R}}
\newcommand{\bQ}{{\mathbb Q}}
\newcommand{\bZ}{{\mathbb Z}}
\newcommand{\bT}{{\mathbb T}}

\newcommand{\cD}{{\mathcal D}}

\newcommand{\cF}{{\mathcal F}}

\newcommand{\cH}{{\mathcal H}}
\newcommand{\cI}{{\mathcal I}}

\newcommand{\cM}{{\mathcal M}}

\newcommand{\cO}{{\mathcal O}}

\newcommand{\cS}{{\mathcal S}}

\newcommand{\Hom}{{\rm{Hom}}}

\newcommand{\wti}{\widetilde}

\newcommand{\Gr}{\text{\rm Gr}}
\newcommand{\DR}{\hbox{\rm DR}}

\newcommand{\lra}{\longrightarrow}

\newcommand{\Int}{{\rm Int}}
\newcommand{\Relint}{{\rm Relint}}
\newcommand{\im}{{\rm Image}}

\newcommand{\ch}{{\rm ch}}

\newcommand{\td}{{\rm td}}

\newcommand{\os}{$O_{\sigma}$}
\newcommand{\vs}{$V_{\sigma}$}

\newcommand{\MHM}{{\rm MHM}}
\newcommand{\MHS}{{\rm MHS}}
\newcommand{\rat}{{\it {rat}}}

\newcommand{\bb}[1]{\mbox{$\mathbb{#1}$}}

\def\sig{\sigma}
\def\Sig{\Sigma}

\def\be{\begin{equation}}
\def\ee{\end{equation}}
\def\bt{\begin{thm}}
\def\et{\end{thm}}
\def\bc{\begin{cor}}
\def\ec{\end{cor}}
\def\br{\begin{rem}}
\def\er{\end{rem}}
\def\bp{\begin{prop}}
\def\ep{\end{prop}}
\def\bl{\begin{lem}}
\def\el{\end{lem}}
\def\bn{\begin{enumerate}}
\def\en{\end{enumerate}}
\def\bex{\begin{example}}
\def\eex{\end{example}}
\def\bd{\begin{df}}
\def\ed{\end{df}}

\begin{document}                        %% Standard LaTeX command

%%      ---------------------------------------------------------------------
%%      -------------------------------- TITLE -----------------------------
%%      ---------------------------------------------------------------------

\title[]{Weighted Ehrhart theory via equivariant toric geometry}
%\dedicatory{Research Announcement}

\author[L. Maxim ]{Lauren\c{t}iu Maxim}
\address{L. Maxim : Department of Mathematics, University of Wisconsin-Madison, 480 Lincoln Drive, Madison WI 53706-1388, USA, \newline
{\text and} \newline Institute of Mathematics of the Romanian Academy, P.O. Box 1-764, 70700 Bucharest, ROMANIA.}
\email {maxim@math.wisc.edu}

\author[J. Sch\"urmann ]{J\"org Sch\"urmann}
\address{J.  Sch\"urmann : Mathematische Institut,
          Universit\"at M\"unster,
          Einsteinstr. 62, 48149 M\"unster,
          Germany.}
\email {jschuerm@math.uni-muenster.de}


\date{\today}

\begin{abstract}
We give a $K$-theoretic and geometric interpretation for a generalized weighted Ehrhart theory of a full-dimensional lattice polytope $P$, depending on a given homogeneous polynomial function $\varphi$ on $P$, and with Laurent polynomial weights $f_Q(y)\in \bZ[y^{\pm 1}]$ associated to the faces $Q \preceq P$ of the polytope. For this purpose, we calculate equivariant $K$-theoretic Hodge--Chern classes of a torus-equivariant mixed Hodge module $\cM$ on the toric variety $X_P$ associated to $P$ (defined via an equivariant embedding of $X_P$ into an ambient smooth variety). 
For any integer $\ell$, we introduce a corresponding equivariant Hodge $\chi_y$-polynomial $\chi_y(X_P, \ell D_P; [\cM])$, with $D_P$ the corresponding ample Cartier divisor  on $X_P$ (defined by the facet presentation of $P$). Motivic properties of the Hodge--Chern classes are used to express this equivariant Hodge polynomial in terms of weighted character sums fitting with a generalized weighted Ehrhart theory. The equivariant Hodge polynomials are shown to satisfy a reciprocity and  purity formula 
fitting with the duality for equivariant mixed Hodge modules, and implying the corresponding properties for the generalized weighted Ehrhart polynomials. In the special case of the equivariant intersection cohomology mixed Hodge module, with the weight function corresponding to Stanley's $g$-function of the polar polytope of $P$, we recover in geometric terms a recent combinatorial formula of Beck--Gunnells--Materov. More generally, motivated by the analogy to the Kazhdan--Lusztig theory, we introduce a duality involution on the free $ \bZ[y^{\pm 1}]$-module of weight functions corresponding to the duality of equivariant mixed Hodge modules, and prove a new reciprocity formula in terms of this duality. This unifies and generalizes the classical reciprocity formula of Brion--Vergne in Ehrhart theory as well as the above-mentioned more recent combinatorial formula of Beck--Gunnells--Materov.
\end{abstract}

\maketitle   

\section{Introduction}
We provide a geometric (equivariant) perspective for the generalized weighted Ehrhart theory developed in \cite{BGM} for a full-dimensional lattice polytope $P$ with a homogeneous polynomial function $\varphi$ defined on it. This is achieved by attaching Laurent polynomial weights $f_Q(y)\in \bZ[y^{\pm 1}]$ of geometric origin to the faces $\emptyset \neq Q \preceq P$ of the polytope (with the partial order given by face inclusion). The Brion--Vergne \cite{BV} combinatorial approach to reciprocity can be linearly extended to this generalized weighted Ehrhart theory. Motivation for using such weights comes from our prior work \cite{CMSSEM} on equivariant Hodge--Chern and Hirzebruch characteristic classes of toric varieties (associated to such polytopes), as well as from the use of Stanley's $g$-functions in the  generalized Ehrhart reciprocity theorem of  Beck--Gunnells--Materov \cite{BGM}.

\subsection{Combinatorial and geometrical posets and their identification} We review here some basic notions explaining the relation between lattice polytopes and toric geometry; for more details see \cite{CLS, Fu2}.

Let $M\cong \bZ^n$ be a lattice and $P \subset M_{\bR}:=M\otimes \bR \cong \bR^n$ be a full-dimensional lattice polytope with vertices in $M$.
To $P$ one associates a fan $\Sig=\Sig_P \subset N_\bR:=N\otimes \bR$, with $N$ the dual lattice of $M$, called the {inner normal fan} of $P$.  To the fan $\Sig_P$ one further associates a projective toric variety $X=X_P$ with torus $\bT\; \simeq (\bC^*)^n$ (whose character lattice is $M$) and an ample Cartier divisor $D=D_P$. There is an order-reversing one-to-one correspondence between the faces $\emptyset \neq Q \preceq P$ of $P$ and the cones $\sig_Q$ of $\Sig_P$, with $\dim_\bR(Q)=n-\dim_\bR(\sig_Q)$. There is an order-reversing one-to-one correspondence between the cones $\sig_Q$ of $\Sig_P$ and the $\bT$-orbits $O_Q$ of the toric variety $X_P$, with $\dim_\bC(O_Q)=n-\dim_\bR(\sig_Q)$. Then the identification between the faces $\emptyset \neq Q \preceq P$ of $P$ and the torus orbits $O_Q$ of the toric variety $X_P$ is order-preserving, with $\dim_\bR(Q)=\dim_\bC(O_Q)$.
The closure $V_Q$ of the orbit $O_Q$ is itself a toric varity for a quotient torus of $\bT$. For $F=Q$ a facet of $P$ (i.e. $\dim_\bR(Q)=n-1$), this gives a $\bT$-invariant prime divisor $D_F:=V_F$.
Since $P$ is full-dimensional, it has a unique facet presentation
$$P=\{m\in M_\bR|\: \langle m, u_F \rangle \geq -a_F \:\text{for all facets $F$ of $P$}\}\:,$$
with $a_F\in \bZ$, $u_F\in N$ the inward pointing facet normal and $ \langle -, - \rangle$  the (real extension of the) duality pairing
$M\times N\to \bZ$. Then by definition $D=D_P:=\sum_F a_F\cdot D_F$ is the  ample Cartier divisor associated to $P$.
Note that the fan $\Sigma_P$ and the toric variety $X_P$ are not changed, if we multiply $P$ by a positive integer $k\in  \bZ_{>0}$,
or translate it by a lattice point $m\in M$:
$$\Sigma_P=\Sigma_{kP-m} \quad \text{and} \quad X_P=X_{kP-m}\:.$$
But the ample Cartier divisor $D_P$ is changed via 
$$D_{kP-m}=k\cdot D_P +div(\chi^m) = k\cdot D_P +\sum_F  \langle m, u_F \rangle \cdot D_F\:,$$
with $div(\chi^m)$ the divisor of the  character $\chi^m$ of $m\in M$ viewed as a rational function on $X_P\supset \bT\to \bC^*$.
For a given face $Q$ of $P$,  choose $m\in Q\cap M$ so that $Q-m$ is a full dimensional lattice polytope in its linear hull
$Span(Q-m)$ relative to the lattice $M'=M\cap Span(Q-m)$. Then the orbit closure $V_Q$ is the toric variety associated to $Q-m$.

If, in addition, $P$ contains the origin in its interior, then its {\it polar polytope} $P^\circ\subset N_\bR$
is defined as in \cite[Section 1.5]{Fu2}, as a full-dimensional rational polytope with respect to the rational points  $N_\bQ\subset N_\bR$, and containing the origin in the interior, given by
$$P^\circ:=\{u\in N_\bR|\: \langle m, u \rangle \geq -1 \:\text{for all $m\in P$}\}\:, \quad \text{with $(P^\circ)^\circ =P$.}$$
 By taking cones at the origin of $N_\bR$ over the proper faces of $P^\circ$, with the empty set $\emptyset$ corresponding to the origin, one gets the same lattice fan $\Sig_P$ (hence the same toric variety $X_P$). This correspondence between the proper faces $\emptyset \preceq Q^\circ \prec P^\circ$ of $P^\circ$ and the cones of $\Sig_P$ is order-preserving, increasing real dimension by $1$. (Here $\dim_\bR \emptyset =-1$ by convention.) This induces an order-reversing one-to-one correspondence between the faces $Q$ of $P$, and  the faces $Q^\circ$ of the polar polytope $P^\circ$, switching the roles of polytopes and emptysets (seen as faces). Moreover,  for a proper face $\emptyset \neq Q \prec P$ of $P$, one has $\dim_\bR(Q)+\dim_\bR(Q^\circ)=n-1$.
In particular, the facets $F$ of $P$ correspond to the vertices $a_F^{-1}\cdot u_F\in N_\bQ$ of $P^\circ$
(with $a_F>0$, since the origin is an interior point of $P$). 

The original polytope $P \subset M_\bR$ is classically used for counting lattice points (and the corresponding Ehrhart theory), fitting with sections of the ample line bundle $\cO_{X_P}(D_P)$ on $X_P$ \cite{Dan,Fu2,CLS}. The torus orbits of $X_P$ give a natural $\bT$-invariant Whitney stratification, which is particularly useful for studying the topology of $X_P$ as well as
in our geometric interpretation of weighted lattice points counting via the theory of Hirzebruch homology classes $T_{y*}$ of \cite{BSY} and the corresponding Hirzebruch--Riemann--Roch theorem; see \cite{MS, MS2}. The polar polytope $P^\circ$ appears in Stanley's work \cite{St} for recursively defining his {\it $g$-polynomials} $g_{Q^\circ}(t) \in \bZ[t]$ (with $g(\emptyset)=1$), for the Eulerian graded poset given by the faces $\emptyset \preceq Q^\circ \preceq P^\circ$,
with rank function $\rho(Q^\circ)=\dim_\bR(Q^\circ)+1$.
 It is well-known \cite{St, DL,F,Sa} (see also \cite{BM}) that these $g$-polynomials are related to the intersection cohomology complex of $X_P$. In the use of these $g$-polynomials, we will (implicitly) switch from $P$ to $P'=kP-m$ for  some $k\in  \bZ_{>0}$ and  
$m\in M$, so that $P'$ contains the origin in its interior. Then $g_{Q^\circ}(t) \in \bZ[t]$ doesn't depend on the choice of $k$ and $m$. 

\subsection{Generalized weighted Ehrhart theory} \label{gwht}
Let $\varphi \colon M_\bR\cong \bR^n \to \bC$ be a homogeneous polynomial function.
We assign to the posets mentioned in the previous subsection Laurent polynomial weights via a {\it weight vector}  
$f=\{f_Q\}$, with $f_Q(y) \in \bZ[y^{\pm 1}]$ indexed here by the non-empty faces $\emptyset \neq Q \preceq P$ of $P$. Motivated by \cite{BGM}, we assemble these weights, the combinatorics of $P$, and the polynomial function $\varphi$ into the {\it generalized weighted Ehrhart polynomial} $E^{\varphi}_{P,f}(\ell, y)$, resp., $\wti{E}^\varphi_{P,f}(\ell, y)$, as follows.
\bd For $\ell \in \bZ_{>0}$, define the generalized weighted Ehrhart polynomial of $P$, $f$ and $\varphi$ by
\be\label{ep1i} 
\begin{split} 
E^\varphi_{P,f}(\ell, y)&:=\sum_{Q \preceq P} f_Q(y) \cdot (1+y)^{\dim(Q)+\deg(\varphi)} \cdot \sum_{m \in  \Relint({{\ell}} Q) \cap M} \varphi(m) \\
&=: (1+y)^{\deg(\varphi)} \cdot \wti{E}^\varphi_{P,f}(\ell, y),
\end{split}
\ee
with  
$ \Relint({{\ell}} Q)$ denoting the relative interior of the face ${\ell} Q$ of the dilated polytope $\ell P$.
\ed

The renormalization $\wti{E}^\varphi_{P,f}(\ell, y)$ is more natural to work with from the geometric perspective considered in this paper.
The Brion--Vergne \cite[Prop.4.1]{BV} combinatorial approach to reciprocity can be linearly extended over $\bZ[y^{\pm 1}]$
to this generalized weighted Ehrhart theory, so that 
$E^\varphi_{P,f}(\ell, y)$, resp., $\wti{E}^\varphi_{P,f}(\ell, y)$, have the following properties:
\begin{enumerate}
\item $E^\varphi_{P,f}(\ell, y)$ is obtained by evaluating a polynomial $E^\varphi_{P,f}(z, y)$ at $z=\ell \in \bZ_{>0}$. Similarly for $\wti{E}^\varphi_{P,f}(\ell, y)$.  In particular, these expressions can be evaluated at any $\ell \in \bZ$.
\item (Reciprocity formula) \ For $\ell \in \bZ_{>0}$, using \cite[Prop.4.1]{BV} for each face $Q$ of $P$, one gets
\be\label{rei} E^\varphi_{P,f}(-\ell, y)=\sum_{Q \preceq P} f_Q(y) \cdot (-1-y)^{\dim(Q)+\deg(\varphi)} \cdot \sum_{m \in  {{\ell}} Q \cap M} \varphi(m).\ee
Similarly, \be\label{recipti} \wti{E}^\varphi_{P,f}(-\ell, y)=\sum_{Q \preceq P} f_Q(y) \cdot (-1-y)^{\dim(Q)} \cdot (-1)^{\deg(\varphi)} \cdot \sum_{m \in  {{\ell}} Q \cap M} \varphi(m).\ee
\item (Constant term) \ For $\ell=0$, \be\label{ctt}E^\varphi_{P,f}(0, y)=\sum_{Q \preceq P} f_Q(y) \cdot (-1-y)^{\dim(Q)+\deg(\varphi)} \cdot \varphi(0),\ee that is, evaluating $\sum_{m \in  {{\ell}} Q \cap M} \varphi(m)$ at $\ell=0$ as $\varphi(0)$ (see \cite[Prop.4.1]{BV}). Similarly,
\be\label{ctt2ti}\wti{E}^\varphi_{P,f}(0, y)=\sum_{Q \preceq P} f_Q(y) \cdot (-1-y)^{\dim(Q)} \cdot (-1)^{\deg(\varphi)} \cdot \varphi(0).\ee
\end{enumerate}

\begin{rem}
If we do not assume that the polynomial function $\varphi$ is homogeneous, then one has for $\ell \in \bZ_{>0}$ the reciprocity formula
$$\wti{E}^\varphi_{P,f}(-\ell, y)=\sum_{Q \preceq P} f_Q(y) \cdot (-1-y)^{\dim(Q)} \cdot \sum_{m \in  {{\ell}} Q \cap M} \varphi(-m)\:,$$
with constant term $\wti{E}^\varphi_{P,f}(0, y)=\sum_{Q \preceq P} f_Q(y) \cdot (-1-y)^{\dim(Q)} \cdot \varphi(0)$.
\end{rem}

In the following, we formulate all our results only for $\varphi$ homogeneous.

 \begin{rem} \label{BVextended} In case $f_Q(y) \in \bZ[y]$ for all $Q$, we are allowed to evaluate these formulae at $y=0$.
For instance, if  $f_Q=\delta_{Q,'Q}$ is a Kronecker function, i.e., $f_Q=1$ for $Q=Q'$ and $f_Q=0$ for $Q\neq Q'$, with $\emptyset \neq Q'$ a given face of $P$, one gets 
$$E^\varphi_{P,f}(\ell, 0) = \wti{E}^\varphi_{P,f}(\ell, 0) = \sum_{m \in  \Relint({{\ell}} Q') \cap M} \varphi(m) \:,$$
with the reciprocity formula (for $\ell \in \bZ_{>0}$)
$$E^\varphi_{P,f}(-\ell, 0) = \wti{E}^\varphi_{P,f}(-\ell, 0) =  (-1)^{\dim(Q') +\deg(\varphi)} \cdot \sum_{m \in  {{\ell}} Q' \cap M} \varphi(m)\:,$$
and with constant term $E^\varphi_{P,f}(0, 0) =\wti{E}^\varphi_{P,f}(0, 0)= (-1)^{\dim(Q') +\deg(\varphi)}  \cdot \varphi(0)$. For $Q'=P$ this is exactly the result of Brion--Vergne
\cite[Prop.4.1]{BV}. But their proof (based on \cite[Prop.3.8]{BV})
can easily be adapted to the case of a general face $Q'$ of $P$ (using a translation $Q'-m$ as a full-dimensional  lattice polytope in $Span(Q'-m)$, together with the projection dual to the inclusion $Span(Q'-m)_\bC\hookrightarrow M_\bC$). Note that $\{\delta_{Q,Q'}|\: \emptyset \neq Q'\preceq P\}$ is a basis of the free $\bZ[y^{\pm 1}]$-module of weight functions.
Finally for $\varphi=1$ and $f_Q=\delta_{Q,Q'}$ we just get the classical Ehrhart polynomial of $\Relint({{\ell}} Q') $:
$$E^1_{P,f}(\ell, 0) = \wti{E}^1_{P,f}(\ell, 0) =  | \Relint({{\ell}} Q') \cap M | \:,$$
with $| - |$ denoting the cardinality of a finite set, together with the reciprocity formula (for $\ell \in \bZ_{>0}$)
$$E^1_{P,f}(-\ell, 0) = \wti{E}^1_{P,f}(-\ell, 0) =  (-1)^{\dim(Q')} \cdot | {{\ell}} Q' \cap M |\:,$$
and constant term $E^1_{P,f}(0, 0) =\wti{E}^1_{P,f}(0, 0)= (-1)^{\dim(Q')}$. In our weighted context it is more natural to start with 
$\Relint({{\ell}} Q') $, whereas classically one starts with the Ehrhart polynomial of $Q'$ (obtained from these formulae by multiplication with 
$(-1)^{\dim(Q')}$).
\end{rem}

Consider now the weight vector given by Stanley's $g$-polynomials \be f_Q(y)=g_{Q^\circ}(-y)=:\wti{g}_Q(-y)\ee for the polar polytope of 
a suitable $P'=kP-m$ (as before), with $g_{\emptyset}(-y)=\wti{g}_P(-y)=1$. In this setup, it was shown in \cite[Thm.1.3, Thm.2.6]{BGM} that the following {\it purity} property holds:
\be\label{pi} E^\varphi_{P,f}(-\ell, y)=(-y)^{\dim(P)+\deg(\varphi)} \cdot E^\varphi_{P,f}(\ell, 1/y).\ee This property can be reformulated as
\be\label{r1bti}
\wti{E}^\varphi_{P,f}(-\ell, y)=(-y)^{\dim(P)} \cdot (-1)^{\deg(\varphi)} \cdot \wti{E}^\varphi_{P, f}(\ell, {1}/{y}).
\ee
Let us note here that our  definition \eqref{ep1i} for $E^\varphi_{P,f}(\ell, y)$ equals $G_{\varphi}(\ell,y)$ from \cite[Formula (14)]{BGM}. 
The constant term $h_{P'^\circ}(y):= E^1_{P,f}(0, -y)= \wti{E}^1_{P,f}(0,- y)$ for $\varphi=1$ is by definition just Stanley's $h$-polynomial of $P'^\circ$, so that one recovers Stanley's ``master duality'' theorem \cite[Thm.3.16.9]{St2}
$$h_{P'^\circ}(y)=y^{\dim(P'^\circ)} \cdot h_{P'^\circ}(1/y) \:.$$
In the geometric terms of the next section, $h_{P'^\circ}(t^2)$ will be the intersection cohomology Poincar\'e polynomial of $X_P$
(as already discussed in \cite{MS2}),
so that this is boiling down to Poincar\'e duality for the intersection cohomology $IH^*(X_P)$ of $X_P$.
In the case when $P$ is a simple polytope, the polar polytope $P'^\circ$ is simplicial, so that $g_{Q^\circ}(-y)=1$, for all faces $Q$ of $P$. In this case, as explained in \cite{BGM}, the equality \eqref{pi} implies the Dehn--Sommerville relations for $P'$, boiling down to Poincar\'e duality for the usual cohomology $H^*(X_P)=IH^*(X_P)$ of $X_P$, which is then a rational homology manifold.\\

Our main new result unifies the {\it reciprocity formulae} of Brion--Vergne
\cite[Prop.4.1]{BV}, extended to all faces $Q'$ of $P$ as in Remark \ref{BVextended}
(which implies \eqref{rei} and \eqref{recipti}), as well as  the {\it purity} results
\eqref{pi} and \eqref{r1bti}, extending the latter also to all faces $Q$ of $P$ (related to the intersection cohomology of the orbit closures
$V_Q$). This will be done in terms of the following duality group homomorphism $\bD$ on the free $\bZ[y^{\pm 1}]$-module of weight functions.

\begin{df}[Duality for weight functions]
The {\it dual} $\bD(f)=\{\bD(f)_Q\}$ of the weight function $f=\{f_Q\}$ is defined as
\be \label{dualweight}
\bD(f)_Q(y):=\sum_{Q \preceq E \preceq P} (1+y)^{\dim_\bR(E)-\dim_\bR(Q)}\cdot (-y)^{-\dim_\bR(E)}\cdot f_E(1/y) \:.
\ee
\end{df}

The geometric meaning of this duality $\bD$ will be explained later on (see formula \eqref{dual2}),
and it is motivated  by the analogy to the {\it Kazhdan--Lusztig theory} (see, e.g., \cite[Sect.4.4]{Pro}, \cite[Thm.7.3.8, Thm.9.3.3]{A2}, and also \cite[Ch.13]{HTT}, \cite{Ta}).
 If we define the standard duality $\bD$ on $\bZ[y^{\pm 1}]$
via $\bD p(y):=p(1/y)$ for $p(y)\in \bZ[y^{\pm 1}]$, then the duality for weight functions has by its definition the {\it module property}
\be
\bD(p\cdot f)= \bD(p)\cdot \bD(f)
\ee
for $p=p(y)\in \bZ[y^{\pm 1}]$ and $f=\{f_Q\}$ a weight function.

\begin{example}[$g$-functions of faces] \label{g-face}
Let $Q'$ be a fixed face of $P$, and choose $k\in  \bZ_{>0}$ and $m\in M$ so that $P'=kQ'-m$ is a full-dimensional lattice polytope in its linear hull $Span(P')$ relative to the lattice $M':=Span(P')\cap M$, containing zero in its relative interior. Let $P'^\circ$ be the polar polytope of $P'\subset M'_\bR$ and define 
$$f_{Q,Q'}(y)= \begin{cases}  g_{Q^\circ}(-y)=:\wti{g}_Q(-y) \quad\text{ for $Q \preceq Q'$,}\\
 \quad\quad 0 \quad \quad \quad \quad \quad \quad \quad\text{otherwise.}
\end{cases}$$
This is independent of the choice of $k$ and $m$, with $f_{Q',Q'}(y)=1$, so that the weight functions $\{f_{Q,Q'}\}$ for all
$\emptyset \neq Q' \preceq P$ are another basis of the free $\bZ[y^{\pm 1}]$-module of weight functions. Moreover
\be\label{duality-g-face}
\bD(f_{Q,Q'}) = (-y)^{-\dim_\bR(Q')}\cdot f_{Q,Q'} \:,
\ee
which is equivalent to \cite[Equation (13)]{BGM} for the polytope $P'$ (as used for $Q'=P$ in the proof of their {\it purity} results).
In particular $\bD$ is an involution.
\end{example}

\begin{thm}[Reciprocity for duality] \label{thm-pur}
Let $P$ be a full-dimensional lattice polytope in $M_\bR$, with $\varphi: M_\bR\simeq \bR^n\to \bC$ a homogenous polynomial
and $f=\{f_Q\}$ a polynomial weight $f_Q(y) \in \bZ[y^{\pm 1}]$ indexed  by the non-empty faces $\emptyset \neq Q \preceq P$ of $P$.
Then one has for $\ell \in \bZ$:
\be\label{rdual} E^\varphi_{P,f}(-\ell, y) = (-y)^{deg(\varphi)}\cdot E^\varphi_{P,\bD(f)}(\ell, 1/y) \:,
\ee
and similarly 
\be\label{rdual2} \wti{E}^\varphi_{P,f}(-\ell, y) = (-1)^{deg(\varphi)}\cdot \wti{E}^\varphi_{P,\bD(f)}(\ell, 1/y) \:.
\ee
\end{thm}

For a geometric statement and proof of the above theorem, see Theorem \ref{repr}. For $f=\{f_{Q,Q'}\}$ as in Example \ref{g-face} we get the corresponding {\it purity} results extending 
\eqref{pi} and \eqref{r1bti} from the case $Q'=P$ to a general face $Q'$ of $P$.
On the other hand, for $f=\{\delta_{Q,Q'}\}$ one gets
$$\bD f_Q(y)= (1+y)^{\dim_\bR(Q')-\dim_\bR(Q)}\cdot (-y)^{-\dim_\bR(Q')}$$
for  $Q \preceq Q'$ and $\bD f_Q(y)=0$ otherwise, with
$$
\wti{E}^\varphi_{P,f}(\ell, y)=(1+y)^{\dim(Q')} \cdot \sum_{m\in \Relint({{\ell}} Q') \cap M} \varphi(m)
$$
and 
$$
\wti{E}^\varphi_{P,\bD(f)}(\ell, \frac{1}{y})=(1+y)^{\dim(Q')} \cdot \left(  (-1)^{\dim(Q')}\cdot  \sum_{m\in {{\ell}} Q'\cap M} \varphi(m) \right).
$$ 
So in this case \eqref{rdual2} is equivalent to the reciprocity formula of Brion--Vergne
\cite[Prop.4.1]{BV}, extended to all faces $Q'$ of $P$ as in Remark \ref{BVextended}.

\br
Note that formulae \eqref{rdual} and \eqref{rdual2} are linear over $\bZ[y^{\pm 1}]$, so it is enough to prove them for a basis of the free $\bZ[y^{\pm 1}]$-module of weight functions. The reciprocity theorem of Brion--Vergne is then equivalent to Theorem \ref{thm-pur}. Similarly, Theorem \ref{thm-pur} is also equivalent to  the purity results  \eqref{pi} and \eqref{r1bti} for the case of a general face of $P$.
\er

As is well known (and recalled in Section \ref{sec6}), the classical geometric origin of Ehrhart theory and the  reciprocity formula of Brion--Vergne
\cite[Prop.4.1]{BV} is the torus-equivariant cohomology of the line bundle $\cO_{X_P}(D_P)$ associated to the ample Cartier divisor $D_P$
on the projective toric variety $X_P$, namely,
\be\label{Intam1}
\Gamma(X_P;\cO_{X_P}(D_P)) = \bigoplus_{m \in P \cap M} \bC \cdot \chi^m \subset \bC[M],
\ee
where $\chi^m$ denotes the \index{character} character defined by $m \in M$, and $ \bC[M]$ is the coordinate ring of the torus $\bT$,
with all higher cohomology trivial by the Demazure vanishing theorem (e.g., see \cite[Prop.4.3.3, Thm.9.2.3]{CLS}).
Moreover, one has by the Batyrev-Borisov vanishing theorem (e.g., see  \cite[Prop.9.2.7]{CLS}) that
\be\label{Intam4}
H^n(X_P;\cO_{X_P}(-D_P)) = \bigoplus_{m \in \Int(P) \cap M} \bC \cdot \chi^{-m} \subset \bC[M],
\ee
with $\Int({P})$ the interior of $P$, and all other cohomology groups vanishing.
Since $X_{\ell P}=X_P$ and $D_{\ell P}=\ell D_P$ for a positive integer $\ell \in \bZ_{>0}$, one gets for the corresponding {\it torus-equivariant Euler characteristics}:
\be\label{Intam5}\begin{split}
\chi^\bT(X_P,\cO_{X_P}(\ell D_P))&=\sum_{m \in \ell P \cap M} \chi^{-m} 
\quad \text{and} \\
\chi^\bT(X_P,\cO_{X_P}(-\ell D_P))&=(-1)^n \cdot \sum_{m \in \Int(\ell P) \cap M} \chi^{m} \:.
\end{split}\ee
Equivariant Serre duality then also gives 
\be\label{Intam6} \begin{split}
\chi^\bT(X_P,\cO_{X_P}(-\ell D_P)  \otimes \omega_{X_P})&=(-1)^n \cdot \sum_{m \in \ell P \cap M} \chi^{m}
\quad \text{and} \\
\chi^\bT(X_P,\cO_{X_P}(\ell D_P) \otimes \omega_{X_P}) &=\sum_{m \in \Int(\ell P) \cap M} \chi^{-m}\:,
\end{split}\ee
with $\omega_{X_P}=\widehat{\Omega}^{\dim(P)}_{X_P}$ the equivariant {dualizing and (Zariski) canonical sheaf} of $X_P$.
Here we choose the convention that  the torus $\bT$ acts on $\bC[M]$ as follows: if $t \in \bT$ and $f \in \bC[M]$, then $t \cdot f \in \bC[M]$ is given by $p \mapsto f(t^{-1} \cdot p)$, for $p \in \bT$ (see \cite[pag.18]{CLS}). In particular, $t \cdot \chi^m = \chi^m(t^{-1}) \chi^m$
(explaing in \eqref{Intam5} the switch from $\chi^m$ to its dual character $\chi^{-m}$).
If we translate  $P$ to $P'=P-m_0$ for a lattice point $m_0\in M$, with $X_{P'}=X_P$, then
$$D_{P'}=D_P+div(\chi^{m_0}), \quad \text{with} \quad  \cO_{X_P}(\ell D_{P'})=\chi^{\ell m_0}\otimes \cO_{X_P}(\ell D_{P}) \:.$$
The multiplication by $\chi^{\ell m_0}$ of the equivariant Euler characteristics on the left hand side above therefore fits on the right hand side with the translation from $P$ to $P'$ via $\chi^{m_0}\cdot \chi^m=\chi^{m_0+m}$.
For a face $Q$ of $P$ and the orbit closure $V_Q\subset X_P$ this implies then similar formulae as above for 
$$\chi^\bT(V_Q,\cO_{X_P}(\pm\ell D_P)|_{V_Q}) \quad \text{ and} \quad
\chi^\bT(V_Q,\cO_{X_P}(\pm\ell D_P)|_{V_Q})\otimes \omega_{V_Q}) \:,$$ 
by just changing $P$ to $Q$, $\Int( - )$ to $\Relint( - )$, and $n$ to $\dim_\bR(Q)$.\\

These {\it  equivariant Euler characteristics} can then be used for the study of the Ehrhart polynomials of
$$ \sum_{m\in {{\ell}} Q \cap M} \varphi(m) \quad  \text{and} \quad \sum_{m\in \Relint({{\ell}} Q) \cap M} \varphi(m)\:,$$
with the reciprocity related to {\it equivariant Serre duality}. Forgetting the additional equivariant structure corresponds to the use of 
$\varphi=1$, giving back the classical Ehrhart polynomials for the number of lattice points 
$$ |{{\ell}} Q \cap M | \quad  \text{and} \quad |\Relint({{\ell}} Q) \cap M |\:.$$

The geometric meaning of the weight functions $f_Q(y)\in \bZ[y^{\pm 1}]$ in our weighted Ehrhart polynomials is of a different nature,
related to the topology of the toric variety $X_P$, e.g., the intersection cohomology $IH^*(X_P)$ of $X_P$,
with the reciprocity related to Poincar\'e duality for the intersection cohomology. The theory of {\it (torus-equivariant) mixed Hodge modules} on the toric variety $X_P$ gives a very natural connection between these two points of view. Their underlying perverse sheaves are related to intersection cohomology and the topology of the toric variety $X_P$, whereas the graded pieces of the filtered de Rham complex of their
underlying filtered D-modules give interesting complexes of coherent $\cO$-sheaves on $X_P$.
Moreover the {\it duality} of  mixed Hodge modules fits with the Verdier duality of perverse sheaves (like intersection cohomology complexes), as well with the  Grothendieck- and (in  the toric context also) Serre-duality of coherent $\cO$-sheaves.
Finally, the stalk or global cohomology (with compact support) of a mixed Hodge module carry a mixed Hodge structure.
Here we only take into account their Hodge filtrations, with the corresponding $\chi_y$-genus (or Hodge-filtered Euler characteristic) of the stalks giving
our weight functions $f_Q(y)\in \bZ[y^{\pm 1}]$. In \cite{MS2} we already studied the weighted Ehrhart polynomials for the special case 
$\varphi=1$ via the theory of mixed Hodge modules on the toric variety $X_P$ in the non-equivariant context,
using their {\it Hirzebruch classes} in the homology of $X_P$. Here we now develop the corresponding 
much reacher torus-equivariant context, which allows the study of the more general weighted Ehrhart polynomials with respect to any homogenous polynomial $\varphi$. Here it is more natural to work with the K-theoretic {\it equivariant Chern classes} of the Hodge modules. Finally, also the explicit {\it duality for weight functions} as discussed before, is new.
As we explain next, it naturally fits with the duality of mixed Hodge modules on the toric variety $X_P$.

\subsection{Geometric viewpoint on generalized weighted Ehrhart theory}
We give  a geometric interpretation and proofs of the  results mentioned before for any homogeneous polynomial function $\varphi$ and any weight vector $f$ via the theory of {\it equivariant} mixed Hodge modules on a projective toric variety and their $K$-theoretic equivariant characteristic classes.  Moreover, using the equivariant intersection Hodge module ${IC}^H_{X_P}$ and its purity, we give a geometric proof of formula \eqref{pi} and its counterpart for orbit closures $V_Q\subset X:=X_P$.

\subsubsection{Geometric realization of weight vectors and duality}
As already indicated in \cite{MS2}, we identify the weight vectors $\{f_Q\}$, seen as $\bZ[y^{\pm 1}]$-valued functions on the set of faces $Q$ of $P$, with the $\bT$-invariant $\bZ[y^{\pm 1}]$-valued constructible functions on $X_P$, via $$\{f_Q\} \mapsto \sum_{Q \preceq P} 1_{O_Q} \cdot f_Q \in F^{\bT}(X)[y^{\pm 1}].$$
Note that $F^{\bT}(X)[y^{\pm 1}]$ is a $\bZ[y^{\pm 1}]$-module freely generated by the $\bT$-invariant constructible functions $1_{O_Q}\in F^{\bT}(X) \subset F^{\bT}(X)[y^{\pm 1}]$,  with $1_{O_{Q'}}$ corresponding to the weight function $\{\delta_{Q,Q'}\}$.  In what follows, we adapt the formulations from  \cite{MS2} in motivic and Hodge-theoretic terms to the equivariant context.

There is a tautological (injective) homomorphism
$$\tau^\bT: F^{\bT}(X)[y^{\pm 1}] \lra K^\bT_0(var/X)[\bL^{-1}]$$
to the equivariant Grothendieck group of $\bT$-varieties over $X$, localized at the Lefschetz motive $\bL=[\bC \to pt]$ (with the trivial $\bT$-action on $\bC$), given by $$1_{O_Q} \cdot f_Q(y)  \mapsto [O_Q \hookrightarrow X] \cdot  f_Q(-\bL).$$ Here, $K^\bT_0(var/X)$ is a $K^\bT_0(var/pt)$-module via the identification $X \times pt \simeq X$. 

In order to identify the equivariant intersection cohomology complex and duality we need to go one step further, and work with 
 the abelian category $\MHM^\bT(X)$ of $\bT$-equivariant mixed Hodge modules  on $X=X_P$. For technical reasons, here we need to choose a $\bT$-equivariant closed embedding $i: X \hookrightarrow Z$ into an ambient smooth $\bT$-variety $Z$,
 and work with  the Grothendieck group  $K_0(\MHM^\bT(X)) =\, K_0(\MHM^\bT_X(Z))$ of equivariant mixed Hodge modules on $Z$ with support on $X$ (so that the corresponding Grothendieck calculus is equivariantly available, see \cite{A}. The main results will be independent of the choice of this embedding $ X \hookrightarrow Z$). There is a canonical homomorphism
$$\chi^\bT_{Hdg}: K^\bT_0(var/X)[\bL^{-1}] \to  K_0(\MHM^\bT(X))= \,K_0(\MHM^\bT_X(Z))$$
such that
$$\chi^\bT_{Hdg}([O_Q \hookrightarrow X] \cdot  f_Q(-\bL))=[(j_Q)_!\bQ^H_{O_Q}] \cdot f_Q(-u),$$
with $j_Q:O_Q \hookrightarrow X \hookrightarrow Z$ the orbit inclusion in the ambient smooth $\bT$-variety $Z$, and $\bQ^H_{O_Q}$ denoting the equivariant constant Hodge module on the orbit $O_Q$.  
Here, $$u=[\bQ(-1)]\in K_0(\MHM^\bT(pt)) \simeq K_0(\MHM(pt))\simeq K_0(\MHS^p),$$ and the module structure on $K_0(\MHM^\bT_X(Z))$ is induced from the exterior product for the identification $Z\simeq pt \times Z$. Also, $\MHS^p$ is the abelian category of graded polarizable $\bQ$-mixed Hodge structures, with $\bQ(-1)$ the corresponding Tate Hodge structure of weight two.

The ambient Grothendieck groups $K^\bT_0(var/X)[\bL^{-1}]$ and $K_0(\MHM^\bT_X(Z))$ are only used formally in the background. For our applications to the toric context, we only need to work with the (sub)groups 
\be\label{gh} \{f_Q\} \longleftrightarrow F^{\bT}(X)[y^{\pm 1}] \longleftrightarrow  \tau^\bT(F^{\bT}(X)[y^{\pm 1}]) 
\longleftrightarrow  \chi^\bT_{Hdg}\left(\tau^\bT(F^{\bT}(X)[y^{\pm 1}])\right)\ee
realizing our weight vectors in the language of constructible functions, motivic, and mixed Hodge module contexts, respectively. These identifications are linear over the following ring coefficient isomorphisms, with $-y=\bL=u$:
$$\bZ[y^{\pm 1}]  \longleftrightarrow \bZ[y^{\pm 1}] \longleftrightarrow \bZ[\bL^{\pm 1}] \longleftrightarrow \bZ[u^{\pm 1}],$$
identifying the coordinate basis elements as follows: $$ e_{Q'}:=\{\delta_{Q,Q'}\}  \longleftrightarrow 1_{O_{Q'}}  \longleftrightarrow 
[O_{Q'} \hookrightarrow X]  \longleftrightarrow [(j_{Q'})_!\bQ^H_{O_{Q'}}]\:.$$

By their definitions, all groups homomorphisms in \eqref{gh} are surjective. In our geometric applications, we work in the group 
$\chi^\bT_{Hdg}\left(\tau^\bT(F^{\bT}(X)\right)$, recovering
 the corresponding weight vector $\{f_Q\}$ via $$f_Q(y):=\chi_y(i_{x_Q}^*[\cM])$$ as the Hodge $\chi_y$-polynomial (see \eqref{316} for a definition) of the stalk class $i_{x_Q}^*[\cM]$, with $[\cM]\in K_0(\MHM^\bT_X(Z))$, and $i_{x_Q}:\{x_Q\} \hookrightarrow Z$ the inclusion (in the ambient smooth $\bT$-variety) of a point $x_Q$ chosen in the orbit $O_Q$. Then our aim is to show that many invariants of $[\cM]$ depend only of this choice of a weight function.

Specific to the toric context is the fact that the stabilizers of the torus action are connected, so that, using \cite[Lem.1.2]{Ta}, any class $[\cM]\in K_0(\MHM^\bT_X(Z))$ can be decomposed as (see Corollary \ref{cmain})
\be\label{f4ni}
\begin{split}
 [\cM]
&=\sum_{Q \preceq P}\, \big[(j_{Q})_!\bQ_{O_Q}^H \big] \cdot i_{x_Q}^*[\cM],
\end{split}
\ee
and the stalk classes do not see the $\bT$-action, i.e., 
$$i_{x_Q}^*[\cM]  \in K_0(\MHM^\bT(pt)) \simeq K_0(\MHM(pt))\simeq K_0(\MHS^p).$$

Since the composed map 
$$\bZ[u^{\pm 1}] \lra K_0(\MHS^p) \overset{\chi_y}{\lra} \bZ[y^{\pm 1}]$$
is the isomorphism given by $u\mapsto -y$, we get the injectivity of the composition
$$F^{\bT}(X)[y^{\pm 1}] \overset{\tau^\bT}{\lra} K^\bT_0(var/X)[\bL^{-1}] \overset{\chi^\bT_{Hdg}}{\lra} K_0(\MHM^\bT_X(Z)).$$
In this way we also see that any weight vector $f$ is realized by the stalk class of a certain equivariant mixed Hodge module complex via $f_Q(y):=\chi_y(i_{x_Q}^*[\cM])$.\\

Even more specific to the toric context is that the submodule $\chi^\bT_{Hdg}\left(\tau^\bT(F^{\bT}(X)[y^{\pm 1}])\right) \subset  K_0(\MHM^\bT_X(Z))$ is preserved under the  duality $\bD^\bT_Z$ of torus-equivariant mixed Hodge modules on $Z$ (with support in $X$), since each torus orbit has a $\bT$-invariant affine neighborhood of product type. 
Note that this duality has the following module property: $\bD^\bT_Z(u \cdot -)=u^{-1}\cdot \bD^\bT_Z(-)$.
As we will see, under the identifications above, this corresponds to our duality involution $\bD$ on the $\bZ[y^{\pm1}]$-module of weight functions.

Let $IC^H_{V_{Q'}}\in \MHM^\bT_X(Z)$ be the equivariant intersection cohomology Hodge module of the orbit closure 
$V_{Q'}\subset X$ corresponding to the face $Q'$ of $P$, and  ${IC'}^H_{V_{Q'}}:=IC^H_{V_{Q'}}[-n']$ with $n'=\dim_\bC(V_{Q'})=\dim_\bR(Q')$.
Since ${IC}^H_{V_{Q'}}$ is pure of weight $n'$, we get that 
$$\bD^\bT_Z{IC'}^H_{V_{Q'}}\simeq {IC'}^H_{V_{Q'}}[2n'](n') \:.$$
This implies the following key purity property:
\be\label{dua} \bD^\bT_Z [{IC'}^H_{V_{Q'}}]=u^{-n'} \cdot 
[{IC'}^H_{V_{Q'}}] \in\chi^\bT_{Hdg}\left(\tau^\bT(F^{\bT}(X)[y^{\pm 1}])\right)\subset  K_0(\MHM^\bT_X(Z)) \:.
\ee
As already explained in the non-equivariant context of \cite{MS2} (for $Q'=P$), the weight vector associated to $[{IC'}^H_{V_{Q'}}]$ is the collection $f_{Q,Q'}$  of Stanley's $g$-functions as in Example \ref{g-face}, i.e., $\wti{g}_Q(-y)=\chi_y(i_{x_Q}^*[{IC'}^H_{V_{Q'}}])$ 
for all $Q  \preceq Q'$ (cf. also \cite{BM,F,DL,Sa}).  So if we use the combinatorial proof of the duality property \eqref{duality-g-face}
of these $f_{Q,Q'}$, then this already implies the identification of the dualites $\bD$ and $\bD^\bT_Z$. We will give later  on (see \eqref{dual2}) another direct geometric proof for the  identification of these dualities, so that \eqref{dua} will imply the key duality property \eqref{duality-g-face}.

\subsubsection{Twisted equivariant Hodge--Chern classes and generalized weighted Ehrhart theory}
We can now explain our geometric interpretation of the generalized weighted Ehrhart theory introduced in Subsection \ref{gwht},  in terms of suitable twists (by the line bundle of the ample Cartier divisor $D_P$ of $P$) of the equivariant Hodge--Chern class transformation $$\DR^\bT_{y}:K_0(\MHM^\bT_X(Z)) \lra K^\bT_{0,X}(Z)[y^{\pm 1}] \simeq K^\bT_0(X)[y^{\pm 1}],$$
generalizing the corresponding non-equivariant Hodge--Chern classes of \cite{BSY,Sc}. 
This equivariant Hodge--Chern transformation is defined in terms of the equivariant filtered De Rham complex of $\cM \in \MHM^\bT_X(Z)$, which equivariantly can only be defined on the ambient smooth $\bT$-variety $Z$. Here, $K^\bT_{0,X}(Z) \simeq K^\bT_0(X)$ is the Grothendieck group of $\bT$-equivariant coherent sheaves on ($Z$ with set-theoretic  support in) $X$.

The equivariant Hodge--Chern class transformation $\DR^\bT_{y}$ commutes with cross product $\boxtimes$ for the identification $Z \times pt \simeq Z$,   fitting with the module structure over $K_0(\MHS^p)$ and the ring homomorphism 
$\chi_y :K_0(\MHS^p) \to\bZ[y^{\pm 1}]$.
Moreover, it is functorial for $\bT$-equivariant morphisms $f\colon (Z,X) \to (Z',X')$ of pairs, with $Z, Z'$ smooth complex algebraic $\bT$-varieties, and $f:Z \to Z'$ restricting to a {\it  proper} morphism $f\colon X \to X'$ of closed $\bT$-invariant subvarieties. 
Especially for $f\colon X \to X'$ an isomorphism, this shows that $\DR^\bT_{y}$ doesn't depend on the choice of the embedding
$X\hookrightarrow Z$.
In this paper, the functoriality is used for the closed inclusions $k_Q: V_{Q} \hookrightarrow Z$ of the closure of the $\bT$-orbits $O_Q$  into the ambient smooth $\bT$-variety $Z$, and their compositions with toric resolutions of these orbit closures, as well as for taking degrees of these classes (i.e., pushing down to a point space)  in case $X$ is compact.
 These properties are then used to prove the following basic identity (see \eqref{drtc}):
\be\label{drtci}
\DR^\bT_y(\big[(j_{Q})_!\bQ^H_{O_Q}\big])=(1+y)^{\dim(O_Q)} \cdot (k_Q)_*[\omega_{V_Q}]_\bT,
\ee
with $\omega_{V_Q}$ the $\bT$-equivariant canonical (and dualizing) sheaf of the orbit closure $V_Q$.

As explained in \cite[Sect.3]{DM}, the equivariant Hodge--Chern class transformation $\DR^\bT_{y}$ commutes with equivariant duality for Grothendieck classes of equivariant mixed Hodge modules on $Z$ with support in $X$, which are in the image of $\chi^\bT_{Hdg}$
 (or, most important for us, contained in $\chi^\bT_{Hdg}\left(\tau^\bT(F^{\bT}(X)[y^{\pm 1}])\right) $), i.e., on such elements one has
 \be\label{Intdume}\DR^\bT_y\circ \bD^\bT_Z=\bD^\bT \circ \DR^\bT_y,\ee 
Here, $\bD^\bT_Z$ is induced from the exact  duality functor of equivariant mixed Hodge modules (on $Z$, with support on $X$).
For $\bD^\bT$, one can use either the equivariant Grothendieck duality on the Grothendieck group $K^\bT_{0,X}(Z)$ induced from Grothendieck duality on $Z$ (remembering the support in $X$), or the equivariant Grothendieck duality on the Grothendieck group $K^\bT_{0}(X)$ induced from Grothendieck duality on $X$ (see \cite{LH} for the foundations of equivariant Grothendieck duality). As in the nonequivariant context, $\bD^\bT$
is extended to $K^\bT_0(X)[y^{\pm 1}]$ by $y \mapsto 1/y$.
With the above notations, this duality property of $\DR_y^\bT$ implies our second basic formula (cf. \eqref{drtz}):
\be\label{drtcdi}
\DR^\bT_y(\big[(j_{Q})_*\bQ^H_{O_Q}\big])=(1+y)^{\dim(O_Q)} \cdot (k_Q)_*[\cO_{V_Q}]_\bT \:.
\ee
Equivalently,
\be\label{drtcdi2}\begin{split}
\DR^\bT_y\left(\bD^\bT_Z\big[(j_{Q})_!\bQ^H_{O_Q}\big]\right) 
&=(-y)^{-\dim(O_Q)}\cdot  (1+y)^{\dim(O_Q)} \cdot (k_Q)_*[\cO_{V_Q}]_\bT \\
&=(-1)^{\dim(O_Q)}\cdot  (1+1/y)^{\dim(O_Q)} \cdot (k_Q)_*[\cO_{V_Q}]_\bT ,
\end{split}\ee
since 
$$\bD^\bT_Z\left( (j_{Q})_!(\bQ^H_{O_Q}) \right)\simeq (j_{Q})_*(\bQ^H_{O_Q})[2 \dim(O_Q)](\dim(O_Q))\:.$$
This allows us to calculate $\DR^\bT_y\left(\bD^\bT_Z( - )\right)$ without knowing the stalk information or weight function 
of  $\big[(j_{Q})_*\bQ^H_{O_Q}\big]$.\\

Let now $D'$ be a $\bT$-invariant Cartier divisor on $X=X_P$, e.g., $D'=\ell D=\ell D_P$ for $\ell \in \bZ$. Let $[\cM] \in K_0(\MHM^\bT_X(Z))$ be fixed.
Then $$\DR^\bT_y([\cM]) \otimes \cO_X(D'):= \DR^\bT_y([\cM]) \cdot [\cO_X(D')] \in K_0^\bT(X)[y^{\pm 1}],$$ where we extend the $K^0_\bT(X)$-module structure of $K_0^\bT(X)$ induced from the tensor product linearly over $\bZ[y^{\pm 1}]$. 
Here $K^0_\bT(X)$ is the Grothendieck group of $\bT$-equivariant vector bundles (or their coherent locally free sheaves of sections) on $X$. For $X$ compact and
$a_X:X \to pt$ the constant map as above, we define the {\it $\bT$-equivariant Hodge polynomial of $(X,D',[\cM])$} as:
\be\label{newi}
\begin{split}
\chi^\bT_y(X,D';[\cM])&:=(a_X)_*\left(\DR^\bT_y([\cM]) \otimes \cO_X(D')\right) \in K_0^\bT(pt)[y^{\pm 1}]=K^0_\bT(pt)[y^{\pm 1}] \\ &=:\chi^\bT(X;\DR^\bT_y([\cM]) \otimes \cO_X(D')) \in \bZ[M][y^{\pm 1}],
\end{split}
\ee
which only depends on the Grothendieck class  $[\cM]$, as well as on $D'$.
Only for $D'=0$ it just depends on  $[\cM]$.
In the second equality of the above definition, we use the isomorphism $K^0_\bT(pt) \simeq \bZ[M]$ induced by the eigenspace decomposition  of a finite dimensional $\bT$-representation
 by characters $\chi^m \in \bZ[M]$ (for $m \in M$ an element in the character lattice of $\bT$), see \eqref{bun}.\\

With these notations, we can now formulate the second main result of this paper (see Theorems \ref{wlpca}, \ref{thmmin}, \ref{t28}):

\bt
\label{wlpcai}
Let $X=X_P$ be the projective toric variety with ample Cartier divisor $D=D_P$ associated to a full-dimensional lattice polytope $P\subset M_{\bR} \cong \bR^n$, and let $\ell \in \bZ_{>0}$. Let $Z$ be an ambient  $\bT$-manifold containing $X$ as a $\bT$-invariant subvariety.
For any face $Q$ of $P$, let $x_Q\in O_{Q}$ be a point in the orbit $O_Q$ of $X$, with inclusion map $i_{x_Q} : \{x_Q\} \hookrightarrow Z$, and set $f_Q(y)=\chi_y(i^*_{x_Q}[\cM])$, with $[\cM]\in K_0(\MHM^\bT_X(Z))$
the (virtual difference) Grothendieck class
 of a $\bT$-equivariant mixed Hodge module on $Z$ with support in $X$. Then we have:
\begin{itemize}
\item[(a)] 
\be\label{wcai}
\begin{split}
\chi^\bT_y(X,\ell D;[\cM]) &= 
 \sum_{Q \preceq P} f_Q(y) \cdot (1+y)^{\dim(Q)} \cdot \sum_{m \in \Relint({{\ell}} Q) \cap M} \chi^{-m},
\end{split}
\ee
with $\Relint$ denoting the relative interior of a face.
\item[(b)] 
\be\label{wcabi}
\begin{split}
\chi^\bT_y(X,-\ell D;[\cM]) &= 
 \sum_{Q \preceq P} f_Q(y) \cdot (-1-y)^{\dim(Q)} \cdot \sum_{m \in {\ell} Q \cap M} \chi^{m},
\end{split}
\ee
\item[(c)] 
\be\label{wcabi0}
\begin{split}
\chi^\bT_y(X,0 \cdot D;[\cM]) &= 
 \sum_{Q \preceq P} \left( f_Q(y) \cdot (-1-y)^{\dim(Q)} \right) \cdot  \chi^{0}.
\end{split}
\ee
\end{itemize}
\et

Here, for $\ell \in \bZ_{>0}$, $\ell D_P=D_{\ell P}$ is the ample Cartier divisor corresponding to the dilated polytope $\ell P$, with the same toric variety $X_{\ell P}=X_P$.

All the formulae of the above theorem rely on the decomposition formula \eqref{f4ni}, together with well-known character decompositions of
 $H^*(X; \cO(\pm\ell D))$ and $H^*(X; \cO(\pm\ell D) \otimes \omega_X)$, for any integer $\ell$, as in the classical work of Danilov \cite{Dan}
  mentioned  in \eqref{Intam1} and \eqref{Intam4}.  More precisely, we use the correspondinging $\bT$-equivariant Euler characteristics $\chi^{\bT}$ as given in  \eqref{Intam5} and  \eqref{Intam6}, and also the mentioned versions for
$$\chi^\bT(V_Q,\cO_{X_P}(\pm\ell D_P)|_{V_Q}) \quad \text{ and} \quad
\chi^\bT(V_Q,\cO_{X_P}(\pm\ell D_P)|_{V_Q})\otimes \omega_{V_Q}) $$
for $Q$ a face of $P$.
In addition, we also use the two basic identities  \eqref{drtci} and \eqref{drtcdi}.
 
Furthermore, using the duality property \eqref{Intdume} for $\DR^\bT_y$ (as in Lemma \ref{laj}), together with equivariant Serre duality, and combining this with $(b)$, we obtain for $\ell \in \bZ_{>0}$
the following {\it main duality formula} for the equivariant Hodge polynomial:
\be\label{dufi}
\chi^\bT_y(X,\ell D;\bD^\bT_Z[\cM])=\chi^{\bT}_{1/y}(X,-\ell D;[\cM])\vert_{m \mapsto -m}.
\ee
Here, the involution $m \mapsto -m$ of $\bZ[M]$ corresponds to the involution $(-)^\vee$ on $K^0_\bT(pt)$ induced from duality of equivariant vector bundles on a point space, i.e.,  finite dimensional representation spaces of $\bT$.
Together with the purity property \eqref{dua}, the duality formula \eqref{dufi} applied to $\cM=  {IC'}^H_{V_Q}$ for
$Q$ an $n'$-dimensional  face of $P$ yields 
  \be\label{icd}
  \chi^\bT_{y}(X,-\ell D;[ {IC'}^H_{V_Q}])=(-y)^{n'} \cdot  \chi^\bT_{1/y}(X,\ell D;[ {IC'}^H_{V_Q}])\vert_{m \mapsto -m} \:.
  \ee
  
  \medskip
  
Let us finally explain the role of the chosen homogeneous polynomial function $\varphi:M_\bR \cong \bR^n \to \bC$ for proving the properties of the generalized weighted Ehrhart polynomials $E^{\varphi}_{P,f}(\ell, y)$ and, resp., $\wti{E}^{\varphi}_{P,f}(\ell, y)$, stated in Subsection \ref{gwht}.  The homogeneous polynomial $\varphi$ defines group homomorphisms $$\wti{E}^\varphi : K^0_\bT(pt)\simeq \bZ[M] \lra \bC, \ \ \  E^\varphi : K^0_\bT(pt)\simeq \bZ[M] \lra \bC[y]$$ via  $\wti{E}^\varphi(\chi^m):= \varphi(-m)$, resp.,
  ${E}^\varphi(\chi^m):= \varphi(-(1+y) \cdot m)$. 
  Since $\varphi$ is homogeneous, the duality involution $m \mapsto -m$ on the left-hand side corresponds to the multiplication by $(-1)^{\deg(\varphi)}$ on the right-hand side (with $y\mapsto y$).
   Linearly extending $\wti{E}^\varphi, E^\varphi $ over Laurent polynomials in $y$, we then have that
   $$\wti{E}^{\varphi}_{P,f}(\ell, y)=
   \wti{E}^\varphi \big(\chi^\bT_y(X,\ell D;[\cM])\big)$$
   and 
   $${E}^{\varphi}_{P,f}(\ell, y)=
   {E}^\varphi \big(\chi^\bT_y(X,\ell D;[\cM])\big)$$
   with $f_Q(y)=\chi_y(i^*_{x_Q}[\cM])$.
 The duality involution $m \mapsto -m$ and $y\mapsto 1/y$ on the left-hand side of the group homomorphisms \begin{equation}\label{efi} \wti{E}^\varphi : K^0_\bT(pt)[y^{\pm 1}] \simeq \bZ[M][y^{\pm 1}] \lra \bC[y^{\pm 1}] , \ \ \  E^\varphi : K^0_\bT(pt)\simeq \bZ[M] \lra \bC[y^{\pm 1}] \end{equation}
corresponds to first switching $y\to 1/y$ and then multiplying by $(-1)^{\deg(\varphi)}$ or, resp., by $(-y)^{\deg(\varphi)}$ on the right-hand side.

Then the formulae for $\chi_y^\bT(X,\ell D; [\cM])$ stated in this subsection translate via  $E^\varphi $ and $\wti{E}^\varphi$ into the corresponding properties for the generalized 
  weighted Ehrhart polynomials $E^{\varphi}_{P,f}(\ell, y)$ and, resp., $\wti{E}^{\varphi}_{P,f}(\ell, y)$, as stated in Subsection \ref{gwht} (see also Section \ref{sec6}).
  
To explain geometrically the polynomial behavior in $\ell$ of $\wti{E}^{\varphi}_{P,f}(\ell, y)$, we have to go from $K$-theory to equivariant homology via the equivariant Todd class transformation, using in addition equivariant localization at the torus fixed points (as explained in the Appendix).

\subsection{Structure of the paper} The paper is structured as follows. 

In Section \ref{sec2}, we give a general introduction to  the calculus of equivariant mixed Hodge modules, together with a definition of the equivariant $K$-theoretic Hodge--Chern class transformation $\DR^\bT_y$ of (classes of) equivariant mixed Hodge modules. We explain here the basic calculus (functoriality, module property, and duality) of these classes, comparing it also with its motivic counterpart from \cite{AMSS, CMSSEM}. 

Preparing the ground for the toric context, in Section \ref{sec3} we compute the equivariant Hodge--Chern classes in terms of a stratification by finitely many torus orbits, reducing them to their motivic counterparts from \cite{CMSSEM}, see Proposition \ref{p8}.

In Section \ref{sec4}, for a $\bT$-invariant Cartier divisor $D$ on $\bT$-invariant compact algebraic sub\-variety $X$ of a smooth ambient $\bT$-variety $Z$, we define for the class $[\cM]\in K_0(\MHM^\bT_X(Z))$  of an equivariant mixed Hodge -module the corresponding equivariant  Hodge polynomial $\chi_y(X,D;[\cM])$ of the triple $(X,D;[\cM])$.

In Section \ref{sec5}, we specialize to the toric context, and make the characteristic class formulae from the previous sections much more explicit.  We also explain the geometric origin of Stanley's $g$-polynomial from the point of view of mixed Hodge modules, and give a geometric explanation of the combinatorial duality \eqref{dualweight}.
Finally, we prove here the two basic formulae \eqref{drtci} and \eqref{drtcdi}, as well as formula \eqref{wcabi0}.

In Section \ref{sec6}, we prove parts $(a)$ and $(b)$ of Theorem \ref{wlpcai}, as well as a technical result (Lemma \ref{laj}) used to obtain the duality formula \eqref{dufi} in the toric context.  We also give in Theorem \ref{repr} a geometric proof via equivariant mixed Hodge modules of the combinatorial reciprocity of Theorem \ref{thm-pur}.

Finally, in the Appendix we explain in a geometric 
way the polynomial behavior in $\ell$ of the generalized weighted Ehrhart polynomials.

\smallskip

\begin{ack} 
We thank the referee for helpful comments and suggestions.\\
L. Maxim acknowledges support from the project ``Singularities and Applications'' - CF 132/31.07.2023 funded by the European Union - NextGenerationEU - through Romania's National Recovery and Resilience Plan.
J. Sch\"urmann was funded by the Deutsche For\-schungs\-gemeinschaft (DFG, German Research Foundation) Project-ID 427320536 -- SFB 1442, as well as under Germany's Excellence Strategy EXC 2044 390685587, Mathematics M\"unster: Dynamics -- Geometry -- Structure. L. Maxim and  J. Sch\"urmann also thank the Isaac Newton Institute for Mathematical Sciences for the support and hospitality during the program ``Equivariant methods in geometry'' when work on this paper was undertaken.
\end{ack}

%%%%%%%%%%%%%%%%%%%

\section{Equivariant Hodge--Chern classes via mixed Hodge modules}\label{sec2}

Let $X$ be a complex algebraic variety with an algebraic $\bT=(\bC^*)^n$-action (e.g., the toric variety $X_P$ associated to a full-dimensional lattice polytope $P \subset M_\bR\cong \bR^n$), and assume that $i: X \hookrightarrow Z$ is a closed $\bT$-equivariant embedding of $X$ into an ambient complex algebraic manifold $Z$ endowed with an algebraic $\bT$-action (if $X$ is a normal quasi-projective $\bT$-variety, such an embedding into a $\bT$-equivariant smooth quasi-projective variety exists by \cite[Thm.1]{Su}, see also \cite[Thm.5.1.25]{CG} for any linear algebraic group action on such a normal quasi-projective variety).

Consider the abelian category $\MHM(X)$ of M. Saito's (algebraic) mixed Hodge modules on $X$, which is stable under duality
$\bD_X$, cross-products $\boxtimes$, closed embeddings $i_*$  and  pullback $f^*[\dim_f]$ for a {\it smooth morphism} of relative dimension $\dim_f$, and compatible with the underlying perverse sheaves via a {\it faithful} functor $rat$.
Here the shift is defined on the corresponding derived category $D^b\MHM(X)$, which has the usual $6$-functor formalism 
compatible with the underlying bounded constructible sheaf complexes of rational vector spaces via
$$rat: D^b\MHM(X)\to D^b_c(X; \bQ)\:,$$
see \cite[Sect.4]{Sa2}.
Then the abelian category $\MHM^{\bT}(X)$ of $\bT$-equivariant mixed Hodge modules on $X$ is defined as in \cite[Def.5.1]{A}
(see also \cite[Sect.6.2]{A2} for the counterpart of perverse sheaves),
together with an exact forgetful functor
$${\rm For} :\MHM^{\bT}(X) \to \MHM(X) $$
(which for the torus $\bT$ and $X$ {\it smooth} is even fully faithful since $\bT$ is connected, see \cite[Thm. 5.2]{A}).
Let $\MHM^{\bT}_X(Z)\subset \MHM^{\bT}(Z)$ be the abelian subcategory  of $\bT$-equivariant mixed Hodge modules on $Z$ with 
support in $X$, i.e., whose underlying  mixed Hodge module ${\rm For}(-)$ is supported on $X$.
Then $i_*: \MHM^{\bT}(X)\to \MHM^{\bT}_X(Z)$ is an {\it equivalence of categories}, since the corresponding non-equivariant equivalence 
of \cite[(4.2.4) and (4.2.7)]{Sa2} commutes with smooth (shifted) pullback. Therefore $i_*$ induces an isomorphism of
the corresponding Grothendieck groups
\be\label{iso-Gr}
i_*: K_0(\MHM^\bT(X)) \simeq K_0(\MHM^\bT_X(Z)) \:,
\ee
which is even a $K_0(\MHM^\bT(pt))$-module isomorphism, with the module structure defined (as before) via cross-products.

 In \cite{A}, Achar introduced a triangulated category $D_{\rm MHM}^{b,\bT}(Z)$, together with a bounded non-degenerate t-structure whose heart is $\MHM^{\bT}(Z)$.  The actual definition of  $D_{\rm MHM}^{b,\bT}(Z)$ is not needed in this paper, we refer to \cite{A} for more details on its construction. However, we note that $D_{\rm MHM}^{b,\bT}(Z)$ is only defined for a smooth $\bT$-variety $Z$, and it is {\it not} the derived category of $\MHM^{\bT}(Z)$. So one has a canonical isomorphism of Grothendieck groups $$K_0(\MHM^{\bT}(Z)) \simeq K_0(D_{\rm MHM}^{b,\bT}(Z)) \:.$$ The category $D_{\rm MHM}^{b,\bT}(Z)$ has the same $6$-functor formalism as developed for $D^b\MHM(Z)$  in \cite{Sa2},  e.g., for a $\bT$-equivariant morphism $f: Z\to Z'$ of smooth
$\bT$-varieties on has 
\be\label{duality-switch}
f_*\circ \bD^\bT_Z \simeq \bD^\bT_{Z'}\circ f_! \quad \text{and} \quad
f^*\circ \bD^\bT_{Z'} \simeq \bD^\bT_Z\circ f^! 
\ee
for the corresponding equivariant duality functors $\bD^\bT_Z$ and 
$\bD^\bT_{Z'}$, which are involutions up to a biduality isomorphism.
All of these functors, as well as the t-structure,  commute with a corresponding  t-exact forgetful functor ${\rm For}: D_{\rm MHM}^{b,\bT}(Z) \to D^b\MHM(Z)$.

Let $D_{\rm MHM, X}^{b,\bT}(Z)\subset D_{\rm MHM}^{b,\bT}(Z)$ be the full triangulated subcategory
of equivariant elements, whose cohomology objects are supported on $X$ (after forgetting the $\bT$-action), i.e., they belong to 
$\MHM^{\bT}_X(Z)$. By this definition, the t-structure above restricts to a t-structure on $D_{\rm MHM, X}^{b,\bT}(Z)$, 
with heart $\MHM^\bT_X(Z)$ as well as
commuting with 
${\rm For}$, and  inducing a canonical isomorphism of Grothendieck groups
(and this is all what we need)
$$K_0(\MHM^\bT(X)) \simeq K_0(\MHM^\bT_X(Z))\simeq K_0(D_{\rm MHM,X}^{b,\bT}(Z)) \:.$$ 

We next discuss a $\bT$-equivariant version of the Hodge--Chern  
class of $X$, see also \cite{DM}. In the non-equivariant case, this was initially introduced in \cite{BSY}, but see also \cite{Sc} for a definition and study via mixed Hodge modules. 
 
Let $X$ be as before a closed $\bT$-invariant subvariety of a complex algebraic $\bT$-manifold $Z$, and let $\cM\in \MHM^\bT_X(Z)$ be  a $\bT$-equivariant mixed Hodge module on $Z$ with support in $X$, e.g., $\cM=IC_X^H$ the $\bT$-equivariant intersection cohomology module (defined  only for $X$ pure-dimensional,  already as an element in $ \MHM^\bT(X)$, by the equivariant intermediate extension of the constant Hodge module from the regular part of $X$,  first to $X$ and then to the ambient $\bT$-variety $Z$).
For $p\in \bZ$ consider 
 de Rham complex $\DR_Z(\cM)$ of the underlying algebraic left $\cD_Z$-module $\cM$ with its integrable connection $\nabla$:
$$ \DR_Z(\cM)=[\cM \overset{ \nabla}{\lra}  \cM\otimes_{\cO_Z} \Omega^1_Z  \overset{ \nabla}{\lra} \cdots \overset{ \nabla}{\lra} \cM\otimes_{\cO_Z} \Omega^{\dim(Z)}_Z]
$$
with $\cM$ in degree $-\dim(Z)$, filtered  by
$$
F_p \DR_Z(\cM) =[F_p\cM \overset{ \nabla}{\lra} F_{p+1}\cM\otimes_{\cO_Z} \Omega^1_Z \overset{ \nabla}{\lra} \cdots  \overset{ \nabla}{\lra} F_{p+ \dim(Z)}\cM\otimes_{\cO_Z} \Omega^{\dim(Z)}_Z]. 
$$
The complex $Gr^F_p\DR_Z(\cM)$ associated to the above  filtered de Rham complex $\DR_Z(\cM)$ is a bounded complex of $\bT$-equivariant coherent sheaves on $Z$ with support in $X$. (We use the same symbol for both the mixed Hodge module and the underlying filtered (left) $\cD_Z$-module.) Moreover, $Gr^F_p\DR_Z(\cM)\simeq 0$ for all but finitely many integers $p$ and $\cM$ fixed (as can be checked already nonequivariantly).

\br\label{rem1}
It follows from \cite[Lem.3.2.6]{Sa1} that the  complex $Gr^F_p\DR_Z(\cM)$ is a bounded complex of (equivariant) coherent $\cO_X$-modules (not just set-theoretically supported on $X$), i.e., these (equivariant) coherent $\cO_X$-modules are annihilated by the ideal sheaf $\cI_X$ of $X$ in $Z$. \qed
\er

The transformation $Gr^F_p\DR_Z$ induces a corresponding map on Grothendieck groups since morphisms of (equivariant) mixed Hodge modules are strict with respect to the Hodge filtration. Let (as before) $D_{{\rm MHM}, X}^{b,\bT}(Z)$ be the full subcategory of $D_{{\rm MHM}}^{b,\bT}(Z)$ consisting of  equivariant elements  with support in $X$. 
For such an equivariant object $\cM \in D_{{\rm MHM}, X}^{b,\bT}(Z)$ on $Z$ with support in $X$, we perform the above construction for each cohomology object of $\cM$. So one gets induced Grothendieck group transformations
$$Gr^F_p\DR_Z: K_0(D_{{\rm MHM}, X}^{b,\bT}(Z)) \simeq K_0(\MHM^\bT_X(Z)) \lra K^\bT_{0,X}(Z) \simeq K^\bT_0(X),$$
where $K^\bT_{0,X}(Z) \simeq K^\bT_0(X)$ is the Grothendieck group of $\bT$-equivariant coherent sheaves on ($Z$ with set-theoretical support in) $X$. The fact that the transformation $Gr^F_p\DR_Z$ takes values in $K^\bT_0(X)$ follows directly from Remark \ref{rem1}. Alternatively, using the ambient variety $Z$, this can also  be deduced from 
the identification 
$$i_*: K^\bT_{0,X}(Z)\simeq K^\bT_0(X)$$ 
induced by the pushforward $i_*:K^\bT_0(X) \to K^\bT_{0,X}(Z)$ for the closed embedding $i:X \hookrightarrow Z$. Indeed, the right inverse of $i_*$ is $i^*:  K^\bT_{0,X}(Z)  \to K^\bT_0(X)$, which is induced from the derived pullback (i.e., alternating sum of equivariant tor-sheaves), since $i$ is of finite tor-dimension (as $Z$ is smooth). So $i_*$ is injective. Moreover, if $\cF$ is an equivariant coherent sheaf on $Z$ supported on $X$, then $\cF$ is annihilated by a finite power of the ideal sheaf $\cI_X$ of $X$, and 
$$[\cF]=\sum_{i \geq 0} [\cI_X^{i} \cdot \cF/{\cI_X^{i +1}}\cdot \cF] \in i_*K^\bT_{0}(X) \subset K^\bT_{0,X}(Z),$$ proving the surjectivity of $i_*$.

We can now make the following definition.
\bd[Equivariant Hodge--Chern transformation]
The {\it $\bT$-equivariant Hodge--Chern class transformation} of a complex algebraic $\bT$-variety $X$ equivariantly embedded in an ambient complex algebraic $\bT$-manifold $Z$ is defined as:
$$ \DR^\bT_y:K_0(\MHM^\bT_X(Z)) \lra K^\bT_0(X)[y^{\pm 1}] $$
\begin{equation*}
\begin{split} 
\DR^\bT_y([\cM])&:=
\sum_{p}\, \big[Gr^F_{-p}\DR_Z([\cM])\big]_\bT \cdot (-y)^p.
\end{split}
\end{equation*}
\ed

\br
By definition, $\DR^\bT_y([\cM])$ depends only on the Grothendieck class $[\cM]$, and not on the choice of a representative $\cM \in D^{b,\bT}_{{\rm MHM},X}(Z)$ for this class. Forgetting the $\bT$-action, one recovers by definition the Hodge--Chern class transformation $\DR_y$ from \cite{BSY, Sc}. \qed
\er

\bp[Functoriality of $\DR^\bT_y$]\label{p22}
Let $f\colon (Z,X) \to (Z',X')$ be a $\bT$-equivariant morphism of smooth complex algebraic $\bT$-varieties $f:Z \to Z'$ restricting to a proper morphism $f\colon X \to X'$ of closed $\bT$-invariant subvarieties. Then $$\DR^\bT_y \circ f_* = f_* \circ \DR^\bT_y: K_0(\MHM^\bT_X(Z)) \to K^\bT_{0,X'}(Z')[y^{\pm 1}].$$
\ep

\begin{proof}
This can be done exactly as explained in \cite[Section 3]{DM}. First, one reduces by induction on the dimension of $X$ to the case when $f\colon X \to X'$ is projective, 
by using for $X$ the equivariant Chow-Lemma of \cite[Thm.2]{Su}, together with a corresponding (graph-) embedding into a $\bT$-equivariant smooth quasi-projective variety, which then exists by \cite[Thm.1]{Su}. For a $\bT$-equivariant morphism 
$f\colon (Z,X) \to (Z',X')$ with  $f\colon X \to X'$ projective one uses then in addition a corresponding {\it strictness} result of
M. Saito \cite[Thm.2.14]{Sa2} for the underlying (non-equivariant) higher direct image mixed Hodge modules,
proved in  \cite[Thm.2.14]{Sa2} directly for a projective morphism $f\colon X \to X'$ of possible singular varieties.
\end{proof}

\br
Under the identifications $K^\bT_{0,X}(Z)\simeq K^\bT_{0}(X)$ and $K^\bT_{0,X'}(Z')\simeq K^\bT_{0}(X')$, the push\-forward $f_*:K^\bT_{0,X}(Z)\to K^\bT_{0,X'}(Z')$ corresponds to the usual pushforward $f_*:K^\bT_{0}(X) \to K^\bT_{0}(X')$. 
\qed \er

\br
The transformation $\DR^\bT_y$ is (in the following sense) independent of the choice of an $\bT$-equivariant embedding of $X$ in a smooth ambient $\bT$-manifold $Z$. In the context of Proposition \ref{p22}, with $f\colon X \to X'$ an isomorphism, the pushforward $f_*\colon K_0(\MHM^\bT_X(Z)) \to K_0(\MHM^\bT_{X'}(Z'))$ is an isomorphism 
fitting with the identification \eqref{iso-Gr}. In particular, the transformation $\DR^\bT_y$ commutes with this isomorphism and the isomorphism $f_*:K^\bT_{0}(X) \to K^\bT_{0}(X')$, as in the above remark. If we now consider two equivariant embeddings $i:X \hookrightarrow Z$ and $i':X \hookrightarrow Z'$, the above observation can be applied to the diagonal embedding $X\hookrightarrow Z \times Z'$ (with the diagonal $\bT$-action) and the projections on the two factors to show that $\DR^\bT_y$ is (up to isomorphism) independent of the choice of an equivariant embedding of $X$ in a smooth ambient $\bT$-manifold. \qed
\er

\bex[Degree of $\DR_y^\bT$]
Let us illustrate the definition of the equivariant Hodge--Chern transformation when $X=Z=pt$ is a point space. In this case, since $\bT$ is connected,  we get by \cite[Lem.1.2]{Ta} that 
$$\MHM^\bT(pt)=\MHM(pt)=\MHS^{p},$$
where $\MHS^{p}$ is the category of graded-polarizable mixed $\bQ$-Hodge
structures, with switching the increasing $\cD$-module filtration to a decreasing Hodge filtration so that on a point space this identification gives $Gr^F_{-p}\DR_Z=Gr^p_F$.

By definition, for $\cM \in D_{\rm MHM}^{b, \bT}(pt)$ we have that $$[\cM]\in 
K_0(D_{\rm MHM}^{b, \bT}(pt))=K_0(\MHM^\bT(pt))=K_0(\MHS^{p}),$$ so $DR^\bT_{y}([\cM])$
 does not see the $\bT$-action, i.e., 
\be\label{316} 
DR^\bT_{y}([\cM])=\sum_{j,p}\,(-1)^j\dim_{\bC} \Gr_F^p(H^j({\rm For}(\cM))\otimes_\bQ \bC)\,(-y)^p
=:\chi_y([{\rm For}(\cM)])\ee
is the corresponding Hodge $\chi_y$-polynomial, defining a ring homomorphism 
$$\chi_y:K_0(\MHS^p) \lra K_0(pt)[y^{\pm 1}]=\bZ[y^{\pm 1}] \hookrightarrow K^\bT_0(pt)[y^{\pm 1}].$$
E.g., for $n \in \bZ$, 
$\chi_y(\bQ(n))=(-y)^{-n}$, where $\bQ(n)$ denotes the Tate Hodge structure of weight $-2n$ on the vector space $\bQ$. 
Similarly, $$\chi_y(X):=\chi_y([H^\bullet_c(X;\bQ_X)])$$ is the (compactly supported) Hodge polynomial of $X$. 
E.g., $\chi_y((\bC^*)^n)=(-y-1)^n$.

By functoriality, if $X$ is compact, and $a_Z: Z \to pt$ is the constant map to a point, we have for 
$\cM\in D^{b,\bT}_{{\rm MHM},X}(Z)$ and $(a_Z)_*\cM \in D_{\rm MHM}^{b, \bT}(pt)$: 
\be\label{eq2} \chi_y(X;[\cM]):=\chi_y([(a_Z)_*\cM])=\chi_y([H^\bullet(X;{\rm For}(\cM))])
\ee
is the degree of the (equivariant) Hodge--Chern class of   $[{\rm For}(\cM)]$ (resp. $[\cM]$) depending only on  $[{\rm For}(\cM)]$. \qed
\eex

\bd
Let $K_0(\MHS^p)^{\rm Tate}$ be the subring of $K_0(\MHS^p)=K_0(\MHM^\bT(pt))$ generated by $u:=[\bQ(-1)]$ and $u^{-1}:=[\bQ(1)]$, with a surjective ring homomorphism $\bZ[u^{\pm 1}] \to K_0(\MHS^p)^{\rm Tate}$. This is an isomorphism since its composition with $\chi_y$ is the isomorphism $\bZ[u^{\pm 1}] \to \bZ[y^{\pm 1}], \ u \mapsto -y$.
\ed

\br Via the description of $\DR^\bT_y$ over a point space, it follows that the transformation $\DR^\bT_y$ commutes with the cross-product  $\boxtimes$ with a point space for the identification $Z \times pt \simeq Z$, since the equivariant Hodge filtration on the product $\cM \boxtimes H$ (with $\cM \in \MHM^\bT_X(Z)$ and $H \in  \MHS^p$) is the product filtration on the underlying filtered $\cD$-module, e.g., see the discussion in \cite[Sections 1.4 and 1.5]{Ta}. Let us also note that the effect of the Tate twist $(n)$ on $\DR^\bT_y$ 
is (like in the case of $\DR_y$) just multiplication by $(-y)^{-n}$, i.e.,  for $\cM \in D^{b,\bT}_{\rm{MHM},X}(Z)$ with $\cM(n):=\cM \boxtimes \bQ(n)$, one has
\be\label{ttw}\DR^\bT_y([\cM(n)])=(-y)^{-n} \cdot \DR^{\bT}_y([\cM]).\ee \qed
\er

\br\label{bla} Note that for $i:X\hookrightarrow Z$ a singular closed embedding into an ambient smooth $\bT$-variety, Achar's $6$-functor calculus for equivariant mixed Hodge modules cannot be used to define the constant mixed Hodge complex on $X$ as a $\bT$-equivariant object, e.g., via the pullback of $\bQ^H_Z$. Nevertheless, its Grothendieck class $[\bQ^H_X]\in K_0(\MHM^\bT_X(Z))$ can be introduced via stratifications as follows. 
Let $\cS_X$ be a $\bT$-invariant algebraic stratification of $X \subset Z$ with smooth strata, and denote by $j_{S,Z}:S \hookrightarrow Z$ the inclusion maps of strata $S \in \cS_X$. Then \be\label{const}[\bQ_X^H]:=\sum_{S \in \cS_X} [(j_{S,Z})_! (j_{S,Z})^*\bQ^H_Z]=\sum_{S \in \cS_X} [(j_{S,Z})_! \bQ^H_S] \in K_0(\MHM^\bT_X(Z)),\ee
where $\bQ^H_S \in D_{\rm MHM}^{b,\bT}(S)$ is the constant mixed Hodge complex on the stratum $S$.
Then, by additivity,  $[\bQ_X^H]$ is independent of the choice of the stratification. Moreover, 
$\chi_y(X;[\bQ^H_X])=\chi_y(X)$ is the Hodge polynomial of $X$, as defined above. Similarly, we can introduce $$i^*:K_0(\MHM^\bT(Z)) \to K_0(\MHM^\bT_X(Z)),$$
with $i^* \circ i_*=id$ on $K_0(\MHM^\bT_X(Z))$ and $i_*:K_0(\MHM^\bT_X(Z)) \to K_0(\MHM^\bT(Z))$ defined by forgetting the support. So $i_*$ is injective. 

Alternatively, if $j:U=Z \setminus X \hookrightarrow Z$ is the inclusion map, one may define the equivariant constant Hodge module complex of $X$ in the ambient $\bT$-manifold $Z$ (up to non-canonical isomorphism) as the cone of the adjunction map
\be\label{q}
\bQ^H_X:=Cone ( j_!  j^*\bQ^H_Z \lra \bQ^H_Z ) \in D^{b,\bT}_{{\rm \MHM}, X}(Z).
\ee
By construction, the equivariant Grothendieck class of $\bQ^H_X$ agrees with \eqref{const}.
\qed
\er

Let us now relate the transformation $\DR^\bT_y$ to its motivic counterpart, as used e.g., in \cite{CMSSEM}. 
Let $X$ be a $\bT$-equivariant complex algebraic variety. The relative equivariant motivic Grothendieck group $K^\bT_0(var/X)$ 
 of varieties over $X$ is the free abelian group on isomorphism classes of $\bT$-varieties over $X$, modulo the usual scissor relations 
 $$[Y \to X]=[U \to X]+[Y\setminus U \to X],$$
 for $U \subset X$ an open $\bT$-invariant subvariety.  If $X=pt$ is a point space, then $K_0^\bT(var/pt)$ is a ring with product given by the cross-product of morphisms for $pt \times X\simeq X$, and the group $K_0^\bT(var/X)$ is a module over $K_0^\bT(var/pt)$  with respect to the cross product. 

For any equivariant morphism $f:Y \to X$ of $\bT$-varieties, there is a well-defined push-forward $f_!:K_0^\bT(var/Y) \to K_0^\bT(var/X)$ defined by composition, which is compatible with the module structure over $K_0^\bT(var/pt)$. \\

Assume now that $i:X\hookrightarrow Z$ is a $\bT$-equivariant closed embedding into an ambient smooth $\bT$-manifold $Z$. Then, using the same argument as in \cite{DM}, there is a natural group homo\-morphism 
$$\chi^\bT_{Hdg}: K^{\bT}_0(var/X) \to K_0(\MHM^\bT_X(Z))$$
uniquely fixed by 
$$[g:Y\to X] \mapsto [(i\circ g)_!\bQ^H_Y],$$
where $Y$ is a smooth $\bT$-variety. For singular $Y$, one decomposes it into smooth $\bT$-invariant subvarieties, and uses the usual scissor relations to define $\chi^\bT_{Hdg}([f:Y\to X])$, e.g., in the notations of Remark \ref{bla}, we have
$$\chi^\bT_{Hdg}([i:X \hookrightarrow Z])=[\bQ^H_X] \in K_0(\MHM^\bT_X(Z)).$$

The transformation $\chi^\bT_{Hdg}$ commutes with pushdown $f_!$ (with proper support) for morphisms of pairs $f: (Z,X) \to (Z',X')$, with $Z, Z'$ smooth. 
Hence, if the restriction $f: X\to X'$ is an isomorphism, it identifies $\chi^\bT_{Hdg}$ for the embedding of $X$ and $X'$, respectively.
Moreover, the composition 
$$mC^\bT_y:=DR^\bT_y \circ \chi^\bT_{Hdg}: K_0^\bT(var/X) \to K_0(\MHM^\bT_X(Z)) \to K_0^\bT(X)[y^{\pm 1}]$$
is independent of the ambient smooth $\bT$-variety $Z$, and it is covariant for $\bT$-equivariant proper morphisms. Moreover, if $X$ is smooth and $\cM=\bQ_X^H \in D_{\rm MHM}^{b,\bT}(X)$ is the constant Hodge module complex on $X$, then the underlying $\bT$-equivariant left $\cD$-module is $\cO_X$ (up to a shift) with the trivial Hodge filtration, so that 
$$Gr^F_{-p}\DR_X(\bQ_X^H)\simeq \Omega^p_X[-p]$$
for $0 \leq p \leq \dim(X)$, and it is $0$ otherwise.  
So, $$mC^\bT_y([id_X])=\sum_{p=0}^{\dim(X)} [\Omega^p_X]_\bT\cdot y^p=:\Lambda_y(\Omega^1_X),$$ the total $\Lambda$-class of the cotangent bundle of $X$.  It  follows that  $mC^\bT_y$ induces the equivariant motivic Chern class transformation introduced in \cite{FRW},
as studied in the toric context also  in \cite{CMSSEM}. For $X$ a quasi-projective $\bT$-variety, this is also 
defined in \cite[Thm.4.2]{AMSS} on the $\bT$-equivariant relative Grothendieck group of $\bT$-equivariant quasi-projective varieties over $X$.\\

In general, if $i:X \hookrightarrow Z$ is a $\bT$-equivariant closed embedding into an ambient smooth $\bT$-variety, we set 
$$mC^\bT_y(X):=mC^\bT_y([id_X])=\DR^\bT_y([\bQ_X^H])$$
for $[\bQ_X^H] \in K^\bT_0(\MHM_X(Z))$ the class defined as in \eqref{const}.
Moreover, if $X$ is pure-dimensional we set
$$IC^\bT_y(X):=\DR^\bT_{y}([{IC'}_X^H]), $$
with  ${IC'}_X^H:=IC_X^H[-\dim(X)]$ the (shifted) equivariant intersection cohomology module of $X$, viewed as a (shifted) $\bT$-equivariant Hodge module on $Z$. 
If $X$ is compact, the degree of $IC^\bT_{y}(X)$ is the $I\chi_y$-polynomial $$I\chi_y(X):=\chi_y(X;{IC'}_X^H).$$ If $X$ is moreover a rational homology manifold (e.g., a complex algebraic $V$-manifold like a simplicial toric variety), then $[{IC'}_X^H] \simeq [\bQ_X^H]$ (which also motivates our shift convention). 

\br\label{prr} As shown in \cite[Thm.4.2]{Sa3}, one has the following identification in $D^b_{\rm coh}(X)$:
\be\label{dub}
\uuline{\Omega}^{p}_X[-p] \simeq Gr_{-p}^F \DR_Z(\bQ^H_X),
\ee
where
$$\uuline{\Omega}^{p}_X:=Gr^p_F(\underline{\Omega}^{\bullet}_X)[p] \in D^b_{\rm coh}(X)$$
are the graded parts of the filtered Du Bois complex $(\underline{\Omega}^{\bullet}_X,F)$  on $X$. Viewing the class $[\bQ^H_X] \in K_0(\MHM^\bT_X(Z))$, and considering $Gr_{-p}^F \DR_Z$ as a transformation on this Grothendieck group, one can take \eqref{dub} as the definition of the equivariant classes $[\uuline{\Omega}^{p}_X]_\bT$ of the graded parts of the Du Bois complex on $X$. This yields the formula 
\be\label{zerob}
mC^\bT_y(X)=\sum_{p = 0}^{\dim(X)} [\uuline{\Omega}_{X}^p]_\bT \cdot y^p \in K^\bT_0(X)[y]
\ee
in the language of equivariant mixed Hodge modules. 
Moreover, for such a $\bT$-invariant closed subvariety $X$ in an ambient smooth $\bT$-variety $Z$, we can use 
the equivariant de Rham complex and Remark \ref{rem1} together with \eqref{q} to define
$$\uuline{\Omega}^{p}_X[-p] \simeq Gr_{-p}^F \DR_Z(\bQ^H_X) \in D^{b,\bT}_{\rm coh}(X),$$
as graded pieces of an equivariant Du Bois complex on $X$.
If $X$ is, moreover, a rational homology manifold, we have ${IC'}^H_X \simeq \bQ^H_X\in D^{b,\bT}_{{\rm \MHM},X}(Z)$, defining an equivariant Du Bois complex on $X$ in terms of the equivariant intersection cohomology module, as used e.g., in \cite{KS}. In the case when $X$ is a complex algebraic $V$-manifold or if $X$ is a  toric variety, $\uuline{\Omega}^{p}_X$ is just the $\bT$-equivariant coherent sheaf (sitting in degree zero, since this is the case nonequivariantly) given by the $\bT$-equivariant Zariski sheaf $\widehat{\Omega}_{X}^p$ of $p$-forms on $X$ (see \cite[p.119 and Prop.4.2]{I}). This holds even equivariantly, since the adjunction map $\uuline{\Omega}^{p}_X \to \widehat{\Omega}_{X}^p=j_*j^*\uuline{\Omega}^{p}_X=j_*\Omega_{X_{\rm reg}}^{p}$ for the open inclusion $j$ of the regular part $X_{\rm reg}$ of $X$ is an equivariant isomorphism, as it is  already nonequivariantly an isomorphism of reflexive sheaves.

In \cite{CMSSEM}, another approach based on simplicial resolutions was used to derive \eqref{zerob}, which in the  toric situation yields (cf. \cite[Prop.3.5]{CMSSEM})
\be\label{zero}
mC^\bT_y(X)=\sum_{p = 0}^{\dim(X)} [\widehat{\Omega}_{X}^p]_\bT \cdot y^p \in K_0(X)[y] .
\ee
\qed
\er

Denote by $\bL:=[\bC \to pt]\in K_0^\bT(var/pt)$ the Lefschetz motive (with $\bT$ acting trivially on $\bC$). Then there is a natural group homomorphism on the localized Grothendieck group
$$\chi^\bT_{Hdg}:K^{\bT}_0(var/X)[\bL^{-1}] \to K_0(\MHM^\bT_X(Z)),$$
mapping $\bL$ to $u=[\bQ(-1)] \in K_0(\MHS^p)=K_0(\MHM^\bT(pt))$. Under the natural forgetful functors, this reduces to the corresponding localized nonequivariant transformation (as used in \cite{MS2}).

\br[Effect of duality]\label{r9}
As shown in \cite[Cor.5.19, Rem.5.20]{Sc}, the non-equivariant transformation $\DR_y$ of \cite{BSY,Sc}
commutes with duality, i.e., \be\label{dum}\DR_y\circ \bD_X=\bD \circ \DR_y,\ee 
where $\bD_X$ is induced from the duality of mixed Hodge modules, and $\bD$ on $K_0(X)$ is induced by Grothendieck duality (e.g., see \cite[Part 1, Sect.7]{FM}, \cite[Part 1]{LH}) extended to $K_0(X)[y^{\pm 1}]$ by $y \mapsto 1/y$. A similar statement holds for the motivic transformation $mC_y$,  
with duality on the localization $K_0(var/X)[\bL^{-1}]$ defined 
as in \cite[Section 4B]{Sc} and \cite{Bit}. All these duality transformations are functorial for proper maps, so that for $X$ compact with $a_X:X \to pt$ proper, one gets
\be\label{eq5}\chi_y(X;\bD_X [\cM])=\chi_{1/y}(X;[\cM]).\ee
Indeed, over a point space, Grothendieck duality and homological duality are trivial, while the mixed Hodge module duality is the usual duality of mixed Hodge structures, e.g., $\bD_{pt}\bQ(n)=\bQ(-n)$.\qed 
\er

As explained in \cite[Sect.3]{DM}, the duality formula \eqref{dum} also holds equivariantly for a (virtual difference) Grothendieck class  $[\cM] \in K_0(\MHM^\bT_X(Z))$ of a $\bT$-equivariant mixed Hodge module on $Z$ with support in $X$ (with $\bD^\bT_Z$ instead of $\bD_X$), at least on (Tate twists of) elements in the image of $\chi^\bT_{Hdg}$ (and this is the only case we'll be interested in later on), namely on such elements one has
 \be\label{dume}\DR^\bT_y\circ \bD^\bT_Z=\bD^\bT \circ \DR^\bT_y.\ee 
Here, $\bD^\bT_Z$ is induced from the exact Verdier duality functor of equivariant mixed Hodge modules (on $Z$, with support on $X$), see \cite{A}. 
For $\bD^\bT$, one can use either the equivariant Grothendieck duality on the Grothendieck group $K^\bT_{0,X}(Z)$ induced from Grothendieck duality on $Z$ (remembering the set-theoretic support in $X$), or the equivariant Grothendieck duality on the Grothendieck group $K^\bT_{0}(X)$ induced from Grothendieck duality on $X$. These are identified by the isomorphism $i_*:K^\bT_{0}(X) \to K^\bT_{0,X}(Z)$, since equivariant Grothendieck duality commutes with proper pushforwards (see \cite{LH} for more details). As in the non\-equivariant context, $\bD^\bT$
is extended to $K^\bT_0(X)[y^{\pm 1}]$ by $y \mapsto 1/y$.

\br
To the authors' knowledge, a corresponding duality transformation on the localized equivariant group $K^{\bT}_0(var/X)[\bL^{-1}]$ is not currently available. But on the image $\tau^\bT(F^{\bT}(X)[y^{\pm 1}]) \subset K^{\bT}_0(var/X)[\bL^{-1}]$ of our group of weight functions it is now available via our duality $\bD$ of weight functions given in the Introduction.
 \qed
\er

\section{Computations via stratifications}\label{sec3}

Let us now indicate how to compute the equivariant Hodge--Chern classes in terms of a $\bT$-invariant stratification given by finitely many $\bT$-orbits.

Let as before $X$ be a closed $\bT$-invariant algebraic subvariety of an algebraic $\bT$-manifold $Z$, with inclusion map $i: X \hookrightarrow Z$. Assume that $X$ has only finitely many $\bT$-orbits, with the strata  $S \in \cS_X$ given by these orbits.
Let $\cM \in D^{b,\bT}_{{\rm MHM},X}(Z)$ be an equivariant element on $Z$ with support in $X$.
Then the underlying bounded constructible complex $$K:=\rat({\rm For}(\cM))\in D^b_c(Z; \bQ)$$ is supported on $X$ and constructible with respect to this orbit stratification, i.e.,
 the sheaves $L_{S,\ell}:=\cH^{\ell}K|_S$ are local systems on $S$ for any $S$ and $\ell \in \bZ$ by the equivariant structure
(see, e.g., \cite[Proof of Prop.6.2.13]{A2}).
Let $U=Z\setminus X$ and $\cS_Z=\cS_X \cup \{U\}$. Then $\cS_Z$ is an algebraic stratification of $Z$ satisfying the above properties for each $S \in \cS_X$, together with $L_{U,\ell}:=\cH^{\ell}K|_U\simeq 0$.
For each $S \in \cS_X$, consider the inclusions
$$j_{S,Z}:S\overset{i_{S,{\bar S}}}{\hookrightarrow} {\bar S} \overset{i_{{\bar S},X}}{\hookrightarrow} X  \overset{i}{\hookrightarrow} Z$$ together with 
$i_{{\bar S},Z}=i \circ i_{{\bar S},X},$
so that $$j_{S,Z}=i_{{\bar S},Z} \circ i_{S,{\bar S}}.$$

Under the above notations and assumptions, and using the $6$-functor formalism from \cite{A} we get as in \cite[Prop.5.1.2]{MSS} the following:
\bp
\label{p7n}
\be\label{f1n} [\cM]  
=\sum_{S\neq U,\ell}\,(-1)^{\ell}\cdot \big[(j_{S,Z})_!L_{S,\ell}^H \big]
\in K_0(\MHM^\bT_X(Z)), 
\ee
where $L_{S,\ell}^H=H^{\ell + \dim (S)} (j_{S,Z})^* \cM[-\dim(S)]$ is the shifted equivariant smooth mixed Hodge module on $S$ with $L_{S,\ell}$ the underlying local system of ${\rm For}(L_{S,\ell}^H)$, and the summation is over strata $S \in \cS_X$.
\ep

\begin{rem}
Note that in this case an orbit-stratum $S$ is {\it affine}, since it is another torus of maybe smaller dimension, so that $(j_{S,Z})_!$ is t-exact
for (equivariant) mixed Hodge modules and  their underlying (equivariant) perverse sheaves. So in this context one has for all strata $S \in \cS_X$:
$$(j_{S,Z})_!L_{S,\ell}^H[\dim(S)]  \in \MHM^{\bT}_X(Z) \simeq \MHM^{\bT}(X) \:.$$
Especially $(j_{S,Z})_!\bQ_{S}^H[\dim(S)] \in \MHM^{\bT}_X(Z) \simeq \MHM^{\bT}(X)$ for all strata $S \in \cS_X$.
\end{rem}

In the notations of Proposition \ref{p7n}, we get the following.

\bc\label{cmain}
Assume, moreover, that for any  orbit $S \in \cS_X$, the stabilizer $\bT_s \subset \bT$ at a point $s \in S$ is connected. Then each $L_{S,\ell}^H \simeq \bQ^H_S \boxtimes L_{S,s,\ell}^H$ is a constant variation of mixed Hodge structures on $S$ with stalk the mixed Hodge structure $L_{S,s,\ell}^H$, and we have
\be\label{f4n}
\begin{split}
 [\cM]
&=\sum_{S,\ell}\,(-1)^{\ell} \cdot
\big[(j_{S,Z})_!\bQ_{S}^H \big] \cdot [L^H_{S, s, \ell}]\\
&=\sum_{S}\, \big[(j_{S,Z})_!\bQ_{S}^H \big] \cdot i_s^*[\cM],
\end{split}
\ee
with $i_s: \{s\} \hookrightarrow Z$ the inclusion, for $s\in S$ chosen, and the summation is over the  orbit-strata $S \in \cS_X$.

In particular,
\be\label{f5n}
\begin{split}
 \DR^\bT_y([\cM]) 
& =\sum_{S\in \cS_X}\,
\DR^\bT_y\big[(j_{S,Z})_!\bQ_S^H \big] \cdot \chi_y(i_s^*[\cM]) \\
&=\sum_{S\in \cS_X}\, mC_y^\bT([j_{S,Z}:S \hookrightarrow Z]) \cdot \chi_y(i_s^*[\cM]) \\
&=\sum_{S\in \cS_X}\,(i_{{\bar S},Z})_*mC^\bT_y ([i_{S,{\bar S}}: S \hookrightarrow {\bar S}])  \cdot \chi_y(i_s^*[\cM]).
 \end{split}
 \ee
\ec

\begin{proof}
The fact that each $L^H_{S,\ell}$ is a constant variation of mixed Hodge structures on $S$ follows by \cite[Lem.1.2]{Ta} using the assumption that the stabilizer $\bT_s$ of $s \in S$ is connected. This proves \eqref{f4n}. To show \eqref{f5n}, we apply the equivariant Hodge--Chern class transformation $\DR_y^\bT$ to \eqref{f4n}, using the fact that $K_0(\MHM^\bT_X(Z))$ is a $K_0(\MHM^\bT(pt))\simeq K_0(\MHS^p)$-module, together with the functoriality of the equivariant motivic Chern class transformation $mC_y^\bT$ for the closed embedding $i_{{\bar S},Z}$.
\end{proof}

\br
Formula \eqref{f4n} shows that under the connected $\bT$-stabilizers assumption, the class $[\cM]$ belongs to the (free) $K_0(\MHM^\bT(pt))=K_0(\MHS^p)$-submodule (respectively, $K_0(\MHS^p)^{\rm Tate}$-submodule, if the stalks $i_s^*[\cM]$ are of Tate type) of $K_0(\MHM^\bT_X(Z))$ generated by classes $[(j_{S,Z})_!\bQ_S^H \big]=\chi^\bT_{Hdg}([j_{S,Z}:S \hookrightarrow Z])$, $S \in \cS_X$.  In what follows, we regard
$$\chi_y(i_s^*[\cM])\in \bZ[y^{\pm 1}]$$ as a weight (depending only on $[\cM]$) attached to the stratum $S \in \cS_X$ containing $s$. In the case of a projective toric variety $X=X_P$ associated to a full dimensional lattice polytope $P$ (as in the Introduction),  
the subgroup $$\chi^\bT_{Hdg}\left(\tau^\bT(F^{\bT}(X)[y^{\pm 1}])\right) \subset K^\bT_0(\MHM_X(Z))$$ 
is the $K_0(\MHS^p)^{\rm Tate}$-submodule generated by classes $[\cM]$ whose stalks are of Tate type:
$$ i_s^*[\cM] \in K_0(\MHS^p)^{\rm Tate} \subset K_0(\MHS^p)\quad \text {for all strata $S$.} $$
In particular, the intersection cohomology Hodge modules of the orbit closures $V_Q$ of $X=X_P$  are of Tate type, so that this property
{\it of Tate type} is stable under duality $\bD^\bT_Z$. This is more generally true for $X$ a quasi-projective toric variety
(as discussed later on). \qed
\er

\br
Under the assumption of connected $\bT$-stabilizers, formula \eqref{f5n} reduces the calculation of equivariant Hodge--Chern classes to the motivic calculus of \cite{CMSSEM}, up to assigning the above mentioned Laurent polynomial weights $\chi_y(i_s^*[\cM])$ to each stratum $S$, with $s \in S$ chosen. For instance, the connected stabilizers assumption is automatically satisfied in the toric context, where such stabilizers are  tori.\qed \er

By functoriality, the terms $mC^\bT_y ([S \hookrightarrow {\bar S}])$ 
appearing in formula \eqref{f5n} are computed via resolutions of singularities as follows (see \cite[Thm.3.3]{CMSSEM},  adapting ideas from \cite[Thm.5.1, Rem.5.2]{We} to the equivariant context):

\bp\label{p8}
Let $i_{S,V}:S\hookrightarrow V$ be a smooth partial $\bT$-equivariant compactification of a stratum $S$ so that 
$D:=V\setminus S$ is a $\bT$-invariant simple normal crossing divisor, and $i_{S,{\bar S}}=\pi_V\circ i_{S,V}$ for a proper morphism
$\pi_V:V\to{\bar S}$. Then we have
\be
mC^\bT_y ([S \hookrightarrow {\bar S}]) = 
(\pi_V)_*\big[ \mathcal{O}_V(-D)  \otimes \Lambda_y \Omega_{V}^1(\log(D))\big]_\bT \in K^\bT_0({\bar S})[y].
\ee
\ep

\br\label{r37} Under the notations and assumption of Proposition \ref{p8}, we also have the following formula:
\be
\DR^\bT_y ([(i_{S,\bar{S}})_*\bQ^H_S]) = 
(\pi_V)_*\big[ \Lambda_y \Omega_{V}^1(\log(D))\big]_\bT \in K^\bT_0({\bar S})[y].
\ee
This can deduced from the formula of Proposition \ref{p8} by equivariant Grothendieck duality, but it will not be used in the following sections.
\er

These formulae become more explicit and much simplified in the toric context, as it will be discussed in Section \ref{sec5} below.

\section{Twisted equivariant Hodge--Chern classes and their degrees}
\label{sec4}

Assume now that a torus $\bT$ with character lattice $M$ acts linearly on a finite dimensional complex vector space $W$, with eigenspace \be\label{eigen} W_{\chi^m}:=\{w\in W \mid t\cdot w=\chi^m(t) w \ \text{ for all} \ t \in \bT\}, \ee for $\chi^m$ a character of $\bT$. Then, as in \cite[Prop.1.1.2]{CLS}, one has an {eigenspace} decomposition $$ W= \bigoplus_{m \in M} W_{\chi^m}.$$ This induces an isomorphism \be\label{bun} K^0_\bT(pt) \overset{\simeq}\lra \bb{Z}[M], \ \ [W]_\bT \mapsto \sum_{m \in M} \dim_\bC W_{\chi^m} \cdot \chi^m,\ee 
identifying the duality involution $(-)^\vee$ on the left hand side with the involution $m \mapsto -m$ on characters.
If $X$ is a compact algebraic $\bT$-variety, with a $\bT$-equivariant coherent sheaf $\cF$, then the cohomology spaces $H^i(X;\cF)$ are finite dimensional $\bT$-representations, vanishing for $i$ large enough. The {\it $K$-theoretic Euler characteristic} of $\cF$ is then defined as: \be\label{Kchi}\begin{split} \chi^\bT(X,\cF)&:=(a_X)_*[\cF]_\bT=[H^\bullet(X;\cF)]_\bT:=\sum_i (-1)^i [H^i(X;\cF)]_\bT \in K^0_\bT(pt)\\ &=\sum_{m \in M} \sum_{i=0}^{\dim(X)} (-1)^i \dim_\bC H^i(X;\cF)_{\chi^m} \cdot \chi^m \in \bZ[M], \end{split} \ee
with $a_X:X \to pt$ the constant map to a point. Note that $\chi^\bT(X,\cF)$ depends only on the class $[\cF]\in K_0^\bT(X)$.

Let now $X$ be a closed $\bT$-invariant compact algebraic subvariety of an algebraic $\bT$-manifold $Z$, 
and let $D$ be a $\bT$-invariant Cartier divisor on $X$. Let $[\cM] \in K_0(\MHM^\bT_X(Z))$ be fixed.
Then $$\DR^\bT_y([\cM]) \otimes \cO_X(D):= \DR^\bT_y([\cM]) \cdot [\cO_X(D)] \in K_0^\bT(X)[y^{\pm 1}],$$ where we extend the $K^0_\bT(X)$-module structure of $K_0^\bT(X)$ induced from the tensor product linearly over $\bZ[y^{\pm 1}]$. For $a_X:X \to pt$ the constant map as above, we can then define the {\it $\bT$-equivariant Hodge polynomial of $(X,D,[\cM])$} as:
\be\label{new}
\begin{split}
\chi^\bT_y(X,D;[\cM])&:=(a_X)_*\left(\DR^\bT_y([\cM]) \otimes \cO_X(D)\right) \in K_0^\bT(pt)[y^{\pm 1}]=K^0_\bT(pt)[y^{\pm 1}] \\ &=:\chi^\bT(X;\DR^\bT_y([\cM]) \otimes \cO_X(D)) \in \bZ[M][y^{\pm 1}],
\end{split}
\ee
which only depends on   the Grothendieck class  $[\cM]$, as well as on $D$.
Only for $D=0$ it just depends on  $[\cM]$.
In other words, $\chi^\bT_y(X,D;[\cM])$ is an incarnation (in terms of the characters of $\bT$) of the degree of the Hodge--Chern class of $[\cM]$ twisted by $\cO_X(D)$.

\br
If $D=0$,  
\be\label{zero}
\begin{split}
\chi^\bT_y(X,0;[\cM])&=(a_X)_* \DR^\bT_y([\cM])=\chi_y(X;[\cM]) \in \bZ[y^{\pm 1}] \\
&= \chi_y(X;[\cM]) \cdot \chi^0 \in \bZ[M][y^{\pm 1}],
\end{split}
\ee
 is the usual Hodge polynomial of $[\cM]$ as in \eqref{eq2}. \qed
\er

\bex If $X$ is smooth and $\cM=\bQ_X^H$ is the constant Hodge module complex (viewed as a Hodge module complex on $Z$), then 
$$\chi^\bT_y(X,D;[\bQ_X^H])=\sum _{p=0}^{\dim(X)} \chi^\bT\left(X, \Omega^p_X \otimes \cO_X(D) \right) \cdot y^p,$$
which is exactly $\chi^\bT_y(X,\cO_X(D))$, the equivariant Hirzebruch polynomial of $D$, e.g., see \cite[Sect.3.2]{CMSSEM}.\qed
\eex

\bex
Assume $X$ is a toric variety or a complex algebraic $V$-manifold. Then
$$\chi^\bT_y(X,D;[\bQ_X^H])=\sum _{p=0}^{\dim(X)} \chi^\bT\left(X, \widehat{\Omega}^p_X \otimes \cO_X(D) \right) \cdot y^p,$$
where  $\widehat{\Omega}_{X}^p$ denotes the corresponding sheaf of Zariski differential $p$-forms, is the Hodge polynomial of $D$ considered in \cite[Sect.3.2]{CMSSEM} in the toric context.  
\qed
\eex

\section{Characteristic class formulae in toric geometry}\label{sec5}

In this section, we specialize the formulae of previous sections to the case of toric varieties (e.g., associated to a full-dimensional lattice polytope), which are known to be Whitney stratified by the orbits of the torus action. 

Let us first introduce some notation. Let $X=X_\Sigma$ be a  toric variety associated to the fan $\Sigma \subset N_\bR\cong N \otimes \bR=\bR^n$, with torus $\bT=N\otimes_\bZ \bC^*=(\bC^*)^n$ having $M$,  the dual of the lattice $N\cong \bZ^n$, as the character lattice. The action of the torus $\bT$ on itself extends to an algebraic action of $\bT$ on $X$ 
with finitely many orbits $O_\sig$, one for each cone $\sig \in \Sig$. In fact,  to a $k$-dimensional cone $\sigma \in \Sig$ there corresponds an $(n-k)$-dimensional torus-orbit \os $ \ \cong (\bC^*)^{n-k}$, and the orbit closure \vs \ of \os \  is itself a toric variety and a $\bT$-invariant subvariety of $X$.
Let $N(\sigma)=N/N_\sigma,$
with $N_\sigma$ denoting the sublattice of $N$ spanned by the points in $\sigma \cap N$.  Let $\bT_{N(\sigma)}=N(\sigma) \otimes_\bZ \bC^*$ be the torus associated to  $N(\sigma)$. For each cone $\nu \in \Sigma$ containing $\sigma$ (in short, $\sig \preceq \nu$), let $\overline{\nu}$ be the image cone in $N(\sigma)_\bR$ under the quotient map $N_\bR \to N(\sigma)_\bR$. Then 
\be \label{star}
Star(\sigma)=\{ \overline{\nu} \subseteq N(\sigma)_\bR \mid \sig \preceq \nu \}
\ee 
is a fan in $N(\sigma)_\bR$, with associated toric variety isomorphic to $V_{\sigma}$ (see, e.g., \cite[Prop.3.2.7]{CLS}), and $\bT$ acts on $V_{\sigma}$ via the morphism \be\label{quot}\bT \to \bT_{N(\sigma)}\ee induced by the quotient map $N \to N(\sigma)$.

Before specializing formula \eqref{f5n} to this toric context, we need the following key additional input derived in \cite{CMSSEM}:

\bp\cite[Prop.3.8]{CMSSEM}\label{p51}
Let $X=X_{\Sigma}$ be a  toric variety defined by the fan $\Sigma$. For any cone $\sigma \in \Sigma$, with orbit $O_{\sigma}$ and inclusion $i_{\sigma} : O_{\sigma} \hookrightarrow V_{\sigma}$ in the orbit closure, we have:
\be\label{drt}
mC^\bT_y([i_\sig])=(1+y)^{\dim(O_{\sigma})} \cdot [\omega_{V_{\sigma}}]_\bT,
\ee
where $\omega_{V_{\sigma}}=\widehat{\Omega}^{\dim(O_\sig)}_{V_{\sigma}}$ is the canonical sheaf on the toric variety $V_{\sigma}$, viewed as a $\bT$-equivariant sheaf via the quotient map $\bT\to \bT_{N(\sig)}$ given by \eqref{quot}.
\ep

\begin{proof}
This follows from Proposition \ref{p8} by using a toric resolution $V$ of $V_\sig$, together with the fact that the locally free sheaf $\Omega_{V}^1(\log(D))$  is  $\bT_{N(\sig)}$-equivariantly, hence also $\bT$-equivariantly, a trivial sheaf of  rank $\dim(O_\sig)$, and there is a $\bT_{N(\sig)}$-equivariant, hence also $\bT$-equivariant isomorphism $R(\pi_V)_*\cO_V(-D)\simeq \omega_{V_\sig}$. 
\end{proof}

In the following we assume that $X=X_\Sigma$ is a quasi-projective toric variety. Then \cite{Su} guarantees the existence of some smooth $\bT$-equivariant ambient variety $Z$ with a closed $\bT$-equivariant embedding $i:X \hookrightarrow Z$. Note also that the torus stabilizer at any point in a torus orbit of $X$ is itself a torus, hence {\it connected}.

Together with formula \eqref{f5n}, Proposition \ref{p51} then gives the following:
\bt\label{t28}
Let $X=X_{\Sigma}$ be a quasi-projective toric variety defined by the fan $\Sigma$, with torus $\bT$. 
For each cone $\sigma \in \Sigma$ with orbit $O_{\sigma}$, let $k_{\sigma} : V_{\sigma} \hookrightarrow Z$ be the inclusion of the orbit closure in the smooth ambient $\bT$-manifold $Z$ mentioned above, and fix a point $x_\sigma \in O_{\sigma}$ with inclusion $i_{x_\sig}: \{x_\sig\} \hookrightarrow Z$.  Let  $[\cM]\in K_0(\MHM^\bT_X(Z))$ be a (virtual difference) Grothendieck class  
 of a $\bT$-equivariant mixed Hodge module on $Z$ with support in $X$. 
Then:
\be\label{drt2g}
\DR^\bT_y([\cM])=\sum_{\sigma \in \Sigma} \, \chi_y(i_{x_\sig}^*[\cM]) \cdot (1+y)^{\dim(O_{\sigma})} \cdot (k_\sigma)_* [\omega_{V_{\sigma}}]_\bT.
\ee
In particular, when $X$ is compact, by taking degrees we get
\be\label{degg}
\chi_y(X;[\cM]) =\sum _{\sig \in \Sig}\, \chi_y(i_{x_\sig}^*[\cM]) \cdot  (-1-y)^{\dim(O_{\sig})}.
\ee
\et 

Theorem \ref{t28} yields the following: 
\bc\label{c28}
Let $X=X_{\Sigma}$ be the quasi-projective toric variety defined by a fan $\Sigma$, with a chosen closed $\bT$-equivariant embedding into a smooth ambient $\bT$-variety $Z$. For each cone $\sigma \in \Sigma$ with orbit $O_{\sigma}$, 
fix a point $x_\sigma \in O_{\sigma}$.   
Then we have the following formula:
\be\label{drt2}
\DR^\bT_y([{IC'}_{X}^H])=\sum_{\sigma \in \Sigma} \, \chi_y({IC'}_{X}^H\vert_{x_\sigma}) \cdot (1+y)^{\dim(O_{\sigma})} \cdot  (k_\sig)_* [\omega_{V_{\sigma}}]_\bT.
\ee
Here, ${IC}_{X}^H \in \MHM^\bT_X(Z)$ is regarded as a $\bT$-equivariant  Hodge module on the ambient $\bT$-manifold $Z$, and each class $[\omega_{V_{\sigma}}]_\bT$ is pushed into $K_{0,X}^\bT(Z)\simeq K_{0}^\bT(X)$.
\ec

\br
In particular, if $X_\Sig$ is a simplicial toric variety then $\chi_y({IC'}_{X}^H\vert_{x_\sigma})=\chi_y(\bQ_{X}^H\vert_{x_\sigma})=1$ for all $\sig \in \Sig$, and one recovers in this case a result from \cite{CMSSEM}. Also note that formula \eqref{drt2} for $X$ an affine toric variety also fits with recent results of the paper \cite{KS}. \qed
\er

Let us next recall the identification of $ \chi_y({IC'}_{X}^H\vert_{x_\sigma})$ with the corresponding Stanley $g$-function, as mentioned in the Introduction, see \cite[Thm.6.2]{DL}, \cite[Thm.1.2]{F}, \cite[Thm.1]{Sa}. Here we follow the treatment from \cite{Sa}, formulated in terms of mixed Hodge modules, as used in this paper.
Recall that to a cone $\sig \in \Sig$ one associates a $\bT$-invariant open affine toric subvariety $U_\sig \subset X_\Sig$, which is of global product type $$U_\sig=Z_\sig \times O_\sig,$$ as follows. Choose a splitting $N=N_\sig \oplus {N'}$ of $N$, where $N_\sig=N \cap \langle \sig \rangle$ is the lattice spanned by $\sig$, with a corresponding splitting $\bT=\bT_\sig \times \bT'$ of $\bT$. Here $Z_\sig$ is the $\bT_\sig$-toric variety of $\sig$ in $N_\sig \otimes \bR$, and $\bT'\simeq O_\sig$. Let $s_\sig:Z_\sig \to U_\sig$ be a trivial section (given by crossing with a point in $\bT'$) of the projection $p_\sig:U_\sig \to Z_\sig$, so that $s_\sig^*{IC'}_{X}^H\simeq {IC'}^H_{Z_\sig}$. Then the stalk cohomology 
of ${IC'}^H_{Z_\sig}$ at the unique $\bT_\sig$-fixed point $x_\sig$ in $Z_\sig$ is  calculated in \cite[(1.7.5)]{Sa} as the primitive intersection cohomology of a projective toric variety $Y$ defined as follows: an interior point of $N_\sig$ gives a $1$-parameter group homomorphism $\bC^* \to \bT_\sig$, so that the $\bC^*$-action on $Z_\sig$ has $x_\sig$ as the only attracting fixed point, and $Y$ is the quotient of $Z_\sig \setminus \{x_\sig\}$ by $\bC^*$. This implies that the stalk cohomologies of ${IC'}_{X}^H$ have pure Hodge structures, which inductively are seen to be of Tate type, so that the odd degree stalk cohomology is vanishing. 
It follows that
$$\chi_y({IC'}_{X}^H\vert_{x_\sigma})=g(-y)=g(t^2),$$
with $-y=t^2$, and $g$ the Poincar\'e polynomial of the primitive intersection cohomology of $Y$.
This $g$-polynomial agrees via the recursive formula \cite[(1.7.9)]{Sa} inductively with Stanley's $g$-function of the corresponding polar polytope (as in the Introduction).
Finally, the above-mentioned Tate property of the stalk cohomologies of intersection cohomology complexes implies that, for a toric variety $X_P$ associated to a full-dimensional lattice polytope $P$, one has 
$$[{IC'}^H_{V_{Q}}] \in\chi^\bT_{Hdg}\left(\tau^\bT(F^{\bT}(X)[y^{\pm 1}])\right)\subset  K_0(\MHM^\bT_X(Z))$$
for any face $Q$ of $P$.\\

With the above notations, we next explain the {\it geometric} origin of the duality transformation \eqref{dualweight} for weight functions. For the use of the equivariant $6$-functor formalism, we need to consider a $\bT_\sig$-equivariant embedding $Z_\sig \hookrightarrow Z'$ into a smooth $\bT_\sig$-variety. Taking products with $\bT'$, this gives a $\bT$-equivariant embedding $U_\sig \hookrightarrow Z' \times \bT'=:Z$ together with a trivial section $s':Z'\to Z$ of the projection $p':Z \to Z'$. Let $\cM \in D_{\rm MHM,U_\sig}^{b,\bT}(Z)$ be a $\bT$-equivariant element on $Z$ with support on $U_\sig$. By induction equivalence (e.g., see \cite[Ex.8.8, Ex.8.9]{A}), we get that $\cM \simeq (p')^*\cM'$ for some $\cM' \in D_{\rm MHM,Z_\sig}^{b,\bT'}(Z')$, with $\cM'\simeq (s')^*\cM$.  Since $p'$ is smooth, we have that
\[
\cM \simeq (p')^!\cM' [-2d'](-d') \in D_{\rm MHM,U_\sig}^{b,\bT'}(Z),
\]
for $d':=\dim(\bT')$. Hence
\[
 (s')^!\cM \simeq  \cM'[-2d'](-d') \simeq (s')^*\cM[-2d'](-d')  \in D_{\rm MHM,Z_\sig}^{b,\bT'}(Z').
\]
Let $i_{x_\sig}:\{x_\sig\} \hookrightarrow Z_\sig \hookrightarrow Z'$ be the point inclusion of the attracting fixed point in $Z_\sig$. Then
\[
(i_{x_\sig})^! (s')^!\cM \simeq \left( (i_{x_\sig})^!\cM'\right)[-2d'](-d').
\]
Let $a:Z'\to \{x_\sig\}$ be the constant map. Applying $a_!$ to the adjunction $(i_{x_\sig})_!(i_{x_\sig})^! \to id$ yields a morphism $(i_{x_\sig})^!\cM' \to a_!\cM'$ in $D_{\rm MHM}^{b,\bT'}(x_\sig)$. The underlying nonequivariant morphism is even an  isomorphism by \cite[Thm.2.10.3]{A2}, because the underlying constructible sheaf complex is $\bT_\sig$-equivariant, and therefore $\bC^*$-equivariant for the $1$-parameter subgroup $\bC^* \to \bT_\sig$ as above (with $x_\sig$ the unique attracting fixed point). This implies that
\[
(i_{x_\sig})^! (s')^![\cM]=a_![\cM'(-d')]  \in K_0(\MHM(x_\sig))=K_0(\MHS^p),
\]
for any $[\cM] \in K_0(D_{\rm MHM,U_\sig}^{b,\bT}(Z))=K_0(\MHM_{U_\sig}^{\bT}(Z)).$ 
In the above notations, for any cone $0 \preceq \tau \preceq \sig$ (i.e., for all $\bT$-orbits $O_\tau$ contained in $U_\sig$, or equivalently, with $O_\sig$ contained in the orbit closure $V_\tau$),
let 
\[
f_\tau(y):=\chi_y([\cM\vert_{x_\tau}]).
\]
Applying $\chi_y$ and using the additivity as in \eqref{f4n}, we get
\begin{equation}\label{dual}
\chi_y\left( (i_{x_\sig})^! (s')^![\cM]   \right)=\sum_{0 \preceq \tau \preceq \sig} (-1-y)^{\dim(O_\tau)-\dim(O_\sig)}\cdot (-y)^{\dim(O_\sig)} \cdot f_\tau(y).
\end{equation}
Dualizing this formula, we get
\begin{equation}\label{dual2}
\begin{split}
f_\sig(\bD^\bT_Z[\cM])&:=\chi_y\left( (i_{x_\sig})^*(s')^*[\bD^\bT_Z\cM]   \right) \\ &=\sum_{0 \preceq \tau \preceq \sig} (-1-1/y)^{\dim(O_\tau)-\dim(O_\sig)}\cdot (-y)^{-\dim(O_\sig)} \cdot f_\tau(1/y) \\
&=\sum_{0 \preceq \tau \preceq \sig} (1+y)^{\dim(O_\tau)-\dim(O_\sig)}\cdot (-y)^{-\dim(O_\tau)} \cdot f_\tau(1/y).
\end{split}
\end{equation} 
Formula \eqref{dual2} matches the definition of the weight function from \eqref{dualweight}, since the calculation of (co)stalks in a given point $x_\sig$ can be done on the $\bT$-invariant open affine neighborhood $U_\sig$.\\

For future reference, we also reformulate here \eqref{drt} in the language of equivariant mixed Hodge modules. Let $k_\sig:V_\sig \hookrightarrow Z$ be the orbit closure embedding into the smooth ambient $\bT$-manifold $Z$ constructed in \cite{Su}. Let $j_\sig:=k_\sig \circ i_\sig:O_\sig \hookrightarrow Z$.
\bc\label{cdrt} In the above notations, 
\be\label{drtc}
\DR^\bT_y(\big[(j_{\sigma})_!\bQ^H_{O_{\sigma}}\big])=(1+y)^{\dim(O_{\sigma})} \cdot (k_\sig)_*[\omega_{V_{\sigma}}]_\bT.
\ee
\ec

\begin{proof}
This follows immediately by applying the proper pushforward $(k_\sig)_*$ to both sides of \eqref{drt}, using the functoriality of $mC^\bT_y$ and the definition of $\chi^\bT_{Hdg}$. Indeed,
$$(1+y)^{\dim(O_{\sigma})} \cdot (k_\sig)_*[\omega_{V_{\sigma}}]_\bT=(k_\sig)_*mC^\bT_y([i_\sig])=mC^\bT_y([j_\sig])=\DR_y^\bT(\big[(j_{\sigma})_!\bQ^H_{O_{\sigma}}] ).$$
\end{proof}

For our applications to generalized weighted Ehrhart theory (or Euler-Maclaurin sums)  in the next section, it is also important to compute $\DR^\bT_y(\big[(j_{\sigma})_*\bQ^H_{O_{\sigma}}\big])$, which, as we will see below, is related to $\DR^\bT_y(\big[(j_{\sigma})_!\bQ^H_{O_{\sigma}}\big])$ via duality.
Before computing $\DR^\bT_y(\big[(j_{\sigma})_*\bQ^H_{O_{\sigma}}\big])$, we need the following lemma.

\bl
For $\sig \in \Sig$ with orbit $O_\sig$ and orbit closure $V_\sig$, we have for the equivariant dualizing complex that
\be\label{Seds}
\bD^\bT ( \cO_{V_{\sigma}} ) \simeq   \omega_{V_\sigma}[\dim(O_\sig)].
\ee
In particular, by equivariant biduality, we get
\be\label{Sed}
\bD^\bT[\omega_{V_\sigma}]_\bT= (-1)^{\dim(O_\sig)} \cdot [\cO_{V_{\sigma}}]_\bT.
\ee
\el

\begin{proof}
Let $\wti{V}_{\sig} \overset{f_{\sig}}{\to} V_{\sig}$ be a toric resolution of singularities of $V_{\sig}$. Then $\wti{V}_{\sig}$ is a toric variety obtained by refining the fan of $V_{\sig}$. Let $O_{\sig} \hookrightarrow \wti{V}_{\sig}$ be the natural open inclusion, with complement the simple normal crossing divisor $D_{\sig}$ whose irreducible components correspond to the rays in the fan of $\wti{V}_{\sig}$. Note that $D_{\sig}$ is a $\bT_{N(\sig)}$-invariant (hence a $\bT$-invariant) divisor. 

Furthermore, note that the canonical sheaf $\omega_{\wti{V}_{\sig}}$ on the toric variety $\wti{V}_{\sig}$ is precisely given by $\cO_{\wti{V}_{\sig}}(-D_{\sig})\simeq \omega_{\wti{V}_{\sig}}$ as $\bT_{N(\sig)}$-equivariant, hence also $\bT$-equivariant sheaves, e.g., see \cite[Thm.8.2.3]{CLS}. And since the toric morphism $f_{\sig}:\wti{V}_{\sig} \to V_{\sig}$ is induced by a refinement of the fan of $V_{\sig}$, it follows from \cite[Thm.8.2.15]{CLS} that there is a $\bT_{N(\sig)}$-equivariant, hence also $\bT$-equivariant,  isomorphism:
\be\label{canonical}
{f_{\sig}}_*\omega_{\wti{V}_{\sig}} \simeq \omega_{V_{\sig}}.
\ee
Altogether, we obtain the following sequence of isomorphisms:
\be
\begin{split}
\bD^\bT(\cO_{V_{\sigma}}) &\simeq \bD^\bT\left((Rf_{\sig})_* \big( \cO_{\wti{V}_{\sigma}} \big)\right)\\ &\simeq (Rf_{\sig})_* \big(\bD^\bT ( \cO_{\wti{V}_{\sigma}}) \big) \\ &\simeq (Rf_{\sig})_* \big( \omega_{\wti{V}_{\sig}}[\dim(O_\sig)] \big) \\ &\simeq \omega_{V_{\sig}}[\dim(O_\sig)],
\end{split}
\ee
where the first isomorphism follows since ${V}_{\sig}$ has only rational singularities (as in the proof of \cite[Prop.3.5]{CMSSEM}); the second isomorphism follows since equivariant Grothendieck duality commutes with proper pushforward (cf. \cite[Thm.25.2]{LH});  the third isomorphism follows 
by \be\label{eqdu} \bD^\bT(\cO_{\wti{V}_{\sig}})\simeq \omega_{\wti{V}_{\sig}}[\dim(O_\sig)],\ee for the smooth variety $\wti{V}_{\sig}$, 
cf. \cite[Thm.28.11]{LH}; finally, the last isomorphism is a consequence of (\ref{canonical}) and the vanishing of the higher derived image sheaves of $\omega_{\wti{V}_{\sig}}$, i.e., 
$R^i(f_{\sig})_*\omega_{\wti{V}_{\sig}}=0$, for all $i >0$, 
see \cite[Thm.9.3.12]{CLS}.

By using equivariant biduality, cf. \cite[Prop.31.9, Lem.31.11]{LH}, we get from \eqref{Seds} that
$$ \cO_{V_{\sigma}}\simeq ( \bD^\bT  \omega_{V_\sigma})[-\dim(O_\sig)],$$
which yields \eqref{Sed}.
\end{proof}

In the notations of the above Corollary \ref{cdrt}, we then have the following
\bp\label{p57} For $j_\sig:O_\sig \hookrightarrow Z$ the orbit inclusion into the ambient $\bT$-manifold $Z$, 
\be\label{drtz}
\DR^\bT_y(\big[(j_{\sigma})_*\bQ^H_{O_{\sigma}}\big])=(1+y)^{\dim(O_{\sigma})} \cdot (k_\sig)_*[\cO_{V_{\sigma}}]_\bT.
\ee
\ep

\begin{proof}
Since $O_\sig$ is smooth, $\bQ^H_{O_{\sigma}}$ is a (shifted) $\bT$-equivariant pure Hodge module of weight $\dim(O_{\sigma})$, so $$\bD^\bT_{O_\sig} \bQ^H_{O_{\sigma}} \simeq \bQ^H_{O_{\sigma}}[2 \dim(O_{\sigma}](\dim(O_{\sigma})).$$
Using \eqref{ttw}, the fact that the even shift does not affect the Grothendieck class,  together with \eqref{duality-switch} for the inclusion $j_{\sigma}$,
we get
\be\label{e1} \DR^\bT_y(\big[(j_{\sigma})_*\bQ^H_{O_{\sigma}}\big])=(-y)^{\dim(O_{\sigma})} \cdot \DR^\bT_y(\bD^\bT_Z \big[ (j_{\sigma})_!\bQ^H_{O_{\sigma}}\big]).\ee
Since $(j_{\sigma})_!\bQ^H_{O_{\sigma}} \in \im(\chi^\bT_{Hdg})$, Remark \ref{r9} gives
\be\label{e2} \begin{split}
\DR^\bT_y(\bD^\bT_Z \big[ (j_{\sigma})_!\bQ^H_{O_{\sigma}}\big])&= \left( \bD^\bT \left( \DR^\bT_y(\big[ (j_{\sigma})_!\bQ^H_{O_{\sigma}}\big])  \right)  \right)_{\vert_{y \mapsto 1/y}}\\
&\overset{\eqref{drtc}}{=} \left( \bD^\bT  \left( (1+y)^{\dim(O_{\sigma})} \cdot (k_\sig)_*[\omega_{V_{\sigma}}]_\bT  \right) \right)_{\vert_{y \mapsto 1/y}}\\
&= (1+\frac{1}{y})^{\dim(O_{\sigma})} \cdot (k_\sig)_*\bD^\bT[\omega_{V_\sigma}]_\bT\\
&\overset{\eqref{Sed}}{=}y^{-\dim(O_\sig)} \cdot (1+y)^{\dim(O_\sigma)} \cdot (-1)^{\dim(O_\sig)} \cdot (k_\sig)_*[\cO_{V_{\sigma}}]_\bT\\
&=(-y)^{-\dim(O_\sig)} \cdot (1+y)^{\dim(O_\sigma)}  \cdot (k_\sig)_*[\cO_{V_{\sigma}}]_\bT.
\end{split}
\ee
Here, in the first two equalities we use a linear extension of equivariant Grothedieck duality over $\bZ[y^{\pm 1}]$. 
Formula \eqref{drtz} is now a consequence of  \eqref{e1} and \eqref{e2}.
\end{proof}

\br
Using Remark \ref{r37}, one can give a direct proof of Proposition \ref{p57}, similar to the proof of Proposition \ref{p51}, using the $\bT$-equivariant isomorphism $R(\pi_V)_*\cO_V \simeq O_{V_\sig}$ coming from the fact that the toric variety $V_\sig$ has rational singularities.
\er

\section{Applications to weighted Ehrhart theory}\label{sec6}

Let now $M\cong \bZ^n$ be a lattice and $P \subset M_{\bR} \cong \bR^n$ be a full-dimensional lattice polytope with associated projective toric variety $X=X_P$, inner normal fan $\Sig=\Sig_P$ and ample Cartier divisor $D=D_P$.
One has (e.g., see  \cite[Prop.4.3.3]{CLS})
\be\label{am1}
\Gamma(X_P;\cO_{X_P}(D_P)) = \bigoplus_{m \in P \cap M} \bC \cdot \chi^m \subset \bC[M],
\ee
where $\chi^m$ denotes the \index{character} character defined by $m \in M$, and $ \bC[M]$ is the coordinate ring of the torus $\bT$.
Also, since $D_P$ is ample, Demazure vanishing (e.g., see \cite[Thm.9.2.3]{CLS}) yields that
\be\label{am2}
H^i(X_P;\cO_{X_P}(D_P))=0, \ {\rm for \ all } \ i>0.
\ee
Moreover, one has by the Batyrev-Borisov vanishing theorem (e.g., see  \cite[Prop.9.2.7]{CLS}) that
\be\label{am3}
H^i(X_P;\cO_{X_P}(-D_P))=0, \ {\rm for \ all } \ i\neq n,
\ee
and 
\be\label{am4}
H^n(X_P;\cO_{X_P}(-D_P)) = \bigoplus_{m \in \Int(P) \cap M} \bC \cdot \chi^{-m} ,
\ee
with $\Int({P})$ the interior of $P$. Also, by equivariant toric Serre duality \cite[Thm.9.2.10]{CLS} (fitting with our description of the equivariant dualizing complex in \eqref{Seds}), 
\be\label{am4d}
\bigoplus_{m \in \Int(P) \cap M} \bC \cdot \chi^{m} = \left( H^n(X_P;\cO_{X_P}(-D_P)) \right)^\vee \simeq H^0(X_P;\cO_{X_P}(D_P) \otimes \omega_{X_P}),
\ee
with $\omega_{X_P}=\widehat{\Omega}^{\dim(P)}_{X_P}$ the equivariant {dualizing and canonical sheaf} of $X_P$, together with $H^i(X_P;\cO_{X_P}(D_P) \otimes \omega_{X_P})=0$ for all $i \neq 0$.

\br\label{eSd}
For any $\bT$-invariant Cartier divisor $D$ on an $n$-dimensional complete toric variety $X$ with torus $\bT$, one has the nondegenerate equivariant pairing
$$H^i(X; \cO_{X}(-D)) \times H^{n-i}(X;\cO_{X}(D) \otimes \omega_X) \lra H^{n}(X;\omega_X) \lra \bC$$
induced by the $\bT$-equivariant sheaf isomorphism
$$\cO_{X}(-D) \otimes \cO_{X}(D) \otimes \omega_X \lra \omega_X$$
and the $\bT$-equivariant trace map $ \omega_X[n] \to \omega_{pt}=\bC$. This induces the 
equivariant Serre duality
$$\left( H^i(X;\cO_{X}(-D)) \right)^\vee \simeq H^{n-i}(X;\cO_{X}(D) \otimes \omega_{X}).$$
\qed
\er

\br 
For a cone $\sig_Q \in \Sig_P$ associated to a face $Q$ of $P$, the corresponding orbit closure $V_{\sig_Q}$ can be identified with the toric variety $X_Q$, with corresponding lattice polytope defined as follows (cf. \cite[Prop.3.2.9]{CLS}): translate $P$ by a vertex $m_0$ of $Q$ so that the origin is a vertex of $Q_0:=Q-m_0$; while this translation by $m_0$ does not change $\Sigma_P$ or $X_P$, $Q_0$ is now a full-dimensional polytope in $Span(Q_0)$ relative to the lattice $Span(Q_0) \cap M$, and  $X_Q$ is the associated toric variety.
However, if one wishes to apply the above formulas to the toric variety $X_Q$ associated to a face $Q$ of $P$, one needs to work with the corresponding divisor \be\label{tran} D_{Q_0}={D_{P-m_0}}\vert_{X_Q}=D_P\vert_{X_Q} + div(\chi^{m_0}).\ee
\qed \er

Formulae \eqref{am1}-\eqref{am4} correspond to the weight decomposition of the respective cohomology groups, e.g., 
\be\label{b30}
H^i(X_P;\cO_{X_P}(D_P))=\bigoplus_{m \in M} H^i(X_P;\cO_{X_P}(D_P))_m,
\ee
as in \cite[Sect.9.1]{CLS}.
To relate this weight decomposition with the eigenspace decomposition of \eqref{eigen}, we recall that  the torus $\bT$ acts on $\bC[M]$ as follows: if $t \in \bT$ and $f \in \bC[M]$, then $t \cdot f \in \bC[M]$ is given by $p \mapsto f(t^{-1} \cdot p)$, for $p \in \bT$ (see \cite[pag.18]{CLS}). In particular, $t \cdot \chi^m = \chi^m(t^{-1}) \chi^m$, so that
$$H^i(X;\cO_X(D))_m = H^i(X;\cO_X(D))_{\chi^{-m}}.$$

Let us note that if for a positive integer $\ell \in \bZ_{>0}$ we consider the dilated polytope $$\ell P:=\{\ell \cdot u \ | \ u \in P \},$$ then the toric varieties corresponding to $P$ and, resp., $\ell P$, are the same, but $D_{\ell P}=\ell D_P$. Together with \eqref{am1}-\eqref{am4d}, this gives for such $\ell \in \bZ_{>0}$:
\be\label{am5}
\chi^\bT(X_P,\cO_{X_P}(\ell D_P))=\sum_{m \in \ell P \cap M} \chi^{-m},
\ee
\be\label{am55}
\chi^\bT(X_P,\cO_{X_P}(-\ell D_P))=(-1)^n \cdot \sum_{m \in \Int(\ell P) \cap M} \chi^{m}
\ee
and
\be\label{am6}
\chi^\bT(X_P,\cO_{X_P}(\ell D_P) \otimes \omega_{X_P})=\sum_{m \in \Int(\ell P) \cap M} \chi^{-m}.
\ee

More generally, we have by \cite[Rem.4.9]{CMSSEM}, applied to the polytope $\ell P$, the following identity:
\be\label{f95}
\chi^\bT(X_P, mC^\bT_y(X_P) \otimes \cO_{X_P}(\ell D_P))=\sum_{Q \preceq P} (1+y)^{\dim(Q)} \cdot \sum_{m \in \Relint(\ell Q) \cap M} \chi^{-m} \in \bZ[M]\otimes_\bZ \bZ[y].
\ee
where $ \Relint$ denotes the relative interior, and the first summation is over the faces $Q$ of the lattice polytope $P$.

For further reference, let us also include here the following.
\bl In the above notations, with $\ell \in \bZ_{>0}$,
\be\label{am56}
\chi^\bT(X_P,\cO_{X_P}(-\ell D_P)  \otimes \omega_{X_P})=(-1)^n \cdot \sum_{m \in \ell P \cap M} \chi^{m}
\ee
\el

\begin{proof}
By toric equivariant Serre duality \cite[Thm.9.2.10]{CLS} and Remark \ref{eSd}, together with \eqref{am1} and \eqref{am2},
\be\label{am4db}
\bigoplus_{m \in P \cap M} \bC \cdot \chi^{-m} = \left( H^0(X_P;\cO_{X_P}(D_P)) \right)^\vee \simeq H^n(X_P;\cO_{X_P}(-D_P) \otimes \omega_{X_P}),
\ee
 and $H^i(X_P;\cO_{X_P}(-D_P) \otimes \omega_{X_P})=0$ for all $i \neq n$.
The assertion follows by reverting to the weight decomposition of these cohomology groups, as indicated above.
\end{proof}

Making use of Theorem \ref{t28}, we can now generalize formula \eqref{f95} as follows. 

\bt
\label{wlpca}
Let $X=X_P$ be the projective toric variety with ample Cartier divisor $D=D_P$ associated to a full-dimensional lattice polytope $P\subset M_{\bR} \cong \bR^n$, and let $\ell \in \bZ_{>0}$. Let $Z$ be an ambient  $\bT$-manifold containing $X$ as a $\bT$-invariant subvariety.
Then for any (virtual-difference) Grothendieck class
$[\cM]\in K_0(\MHM^\bT_X(Z))$ of a mixed Hodge module  on $Z$ with support in $X$, we have 
\be\label{wca}
\begin{split}
\chi^\bT_y(X,\ell D;[\cM]) &= 
 \sum_{Q \preceq P} \chi_y(i_{x_Q}^*[\cM]) \cdot (1+y)^{\dim(Q)} \cdot \sum_{m \in \Relint({{\ell}} Q) \cap M} \chi^{-m},
\end{split}
\ee
where the first summation on the right is over the faces $Q$ of $P$, and $x_Q\in O_{\sigma_Q}$ is a point in the orbit of $X$ associated to (the cone $\sigma_Q$ in the inner normal fan of $P$ corresponding to) the face $Q$, with inclusion map $i_{x_Q} : \{x_Q\} \hookrightarrow Z$.
\et

\begin{proof} For a face $Q$ of $P$, denote by $k_Q:V_{\sig_Q}:=X_Q \hookrightarrow X \hookrightarrow Z$ the inclusion in the ambient $\bT$-manifold $Z$ of the orbit closure associated to the (cone $\sigma_Q$ of the) face $Q$. Note that we have $\dim_\bR(Q)=\dim_\bR(O_{\sig_Q}).$ Let $\bT_Q$ be the quotient torus of $\bT$ corresponding to $X_Q$, as in \eqref{quot}.

In what follows, we will denote all constant maps to a point by $a$, regardless of their domains (which should be understood from the context).
Then by \eqref{new} and \eqref{drt2g} the following  equality holds:
\begin{equation*} \begin{split} \chi^\bT_y(X,\ell D;\cM)&=a_*\left(\DR^\bT_y([\cM]) \otimes \cO_X(\ell D)\right) \\
&= \sum_{Q \preceq P} \chi_y(i_{x_Q}^*[\cM]) \cdot  (1+y)^{\dim(Q)} \cdot a_* \left( (k_Q)_*[\omega_{X_Q}]_\bT \otimes \cO_{X}(\ell D) \right).\end{split} \end{equation*}
So it remains to prove that, for any face $Q$ of $P$, we have that:
\be\label{lem2}
\chi^\bT(X,(k_Q)_*[\omega_{X_Q}]_\bT \otimes \cO_{X}(\ell D))=a_* \left( (k_Q)_*[\omega_{X_Q}]_\bT \otimes \cO_{X}(\ell D) \right) = \sum_{m \in \Relint(\ell Q) \cap M} \chi^{-m}.
\ee

This follows from the projection formula and \eqref{am6}, by making use of formula \eqref{tran}. Indeed,
\begin{equation*}
\begin{split}
a_* \left( (k_Q)_*[\omega_{X_Q}]_\bT \otimes \cO_{X}(\ell D) \right) 
 &= a_*  (k_Q)_* \left([\omega_{X_Q}]_\bT \otimes k_Q^*\cO_{X}(\ell D) \right) \\
&\overset{(\ref{tran})}{=} a_*  \left([\omega_{X_Q}]_\bT \otimes \cO_{X_Q}(\ell D_{Q_0}-div(\chi^{\ell m_0})) \right) \\
&=\chi^{-\ell m_0} \cdot a_*\left( [\omega_{X_Q}]_\bT \otimes \cO_{X_Q}(\ell D_{Q_0}) \right) \\
&\overset{({\ast})}{=} 
\chi^{-\ell m_0} \cdot a_*\left( [\omega_{X_Q}]_{\bT_Q} \otimes \cO_{X_Q}(\ell D_{Q_0}) \right) \\
&\overset{(\ref{am6})}{=} 
\chi^{-\ell m_0} \cdot \left(   
\sum_{m \in \Relint \ \ell (Q-m_0) \cap M} \chi^{-m}  \right)\\
&=\sum_{m \in \Relint(\ell Q) \cap M} \chi^{-m},
 \end{split}
\end{equation*}
where $(\ast)$ uses \cite[Prop.3.11]{CMSSEM}, showing  the equality of the equivariant classes 
\begin{equation}\label{equality} [\omega_{X_Q}]_\bT =  [\omega_{X_Q}]_{\bT_Q}\end{equation}
under the natural inclusion of Grothendieck groups $$K_0^{\bT_Q}(X_Q) \hookrightarrow K_0^{\bT}(X_Q)$$ induced by the surjection $\bT \to \bT_Q$. Here $m_0$ is a vertex of $Q$ as in (\ref{tran}), so that $Q_0:=Q-m_0$ is a full-dimensional lattice polytope  in $Span(Q_0)$ relative to the lattice $Span(Q_0)\cap M$, with $X_Q$ the associated toric variety.
\end{proof}

\br
For future use, we include here the following formula, which one gets as in the proof of  \eqref{lem2}, but using $\cO_{X_Q}$ instead of $\omega_{X_Q}$, together with \eqref{am5}:
\be\label{lem2b}
\chi^\bT(X,(k_Q)_*[\cO_{X_Q}]_\bT \otimes \cO_{X}(\ell D))=
a_* \left( (k_Q)_*[\cO_{X_Q}]_\bT \otimes \cO_{X}(\ell D) \right) = \sum_{m \in \ell Q \cap M} \chi^{-m}.
\ee
\qed
\er

\bd[Generalized weighted Ehrhart polynomials]\label{def65}
Let $\varphi:M_\bR \cong \bR^n \to \bC$ be a homogeneous polynomial function, and $\ell \in \bZ_{>0}$. 
With the above notations, the expression
\be\label{whpd}E^\varphi_{P,\cM}(\ell, y):=\sum_{Q \preceq P} \chi_y(i_{x_Q}^*[\cM]) \cdot (1+y)^{\dim(Q)+\deg(\varphi)} \cdot \sum_{m \in  \Relint({{\ell}} Q) \cap M} \varphi(m)\ee
is obtained from $\chi^\bT_y(X,\ell D;\cM)$ by substituting  $\chi^{-m}$ by $\varphi\left((1+y)m\right)=(1+y)^{\deg(\varphi)}\cdot \varphi(m)$, for each $m$.
As we will show later on in the Appendix (but see also \cite[Prop.4.1]{BV}), $E^\varphi_{P,\cM}(\ell, y)$ is obtained by evaluating a polynomial $E^\varphi_{P,\cM}(z, y)$ at $z=\ell$. We call this polynomial the {\it generalized weighted Ehrhart polynomial} of $P$ and $\varphi$ with coefficients induced from $[\cM]\in K_0(\MHM^\bT_X(Z))$. Similarly, one can define \be\label{whpdt}\wti{E}^\varphi_{P,\cM}(\ell, y):=\sum_{Q \preceq P} \chi_y(i_{x_Q}^*[\cM]) \cdot (1+y)^{\dim(Q)} \cdot \sum_{m \in  \Relint({{\ell}} Q) \cap M} \varphi(m)\ee
which obtained from $\chi^\bT_y(X,\ell D;\cM)$ by substituting  $\chi^{-m}$ by $\varphi\left(m\right)$, for each $m$. Since $\varphi$ is homogeneous, this substitution fits with the duality involution $\chi^m \mapsto \chi^{-m}$ for $m \in M$, via $\varphi(-m)=(-1)^{\deg(\varphi)} \cdot \varphi(m)$.

More generally, for a Laurent polynomial {\it weight function} $$f:\{\text{faces of } P\} \to \bZ[y^{\pm 1}], \ \ \ Q \mapsto f_Q(y)$$ defined on the set of faces of $P$, we define as in the Introduction an associated generalized weighted Ehrhart polynomial by
\be\label{whpg} 
\begin{split}
E^\varphi_{P,f}(\ell, y)&:=\sum_{Q \preceq P} f_Q(y) \cdot (1+y)^{\dim(Q)+\deg(\varphi)} \cdot \sum_{m \in  \Relint({{\ell}} Q) \cap M} \varphi(m)\\
&=:(1+y)^{\deg(\varphi)} \cdot \wti{E}^\varphi_{P,f}(\ell, y).\end{split}\ee
\ed

\br\label{r34}
Note that any Grothendieck class $[\cM]\in K_0(\MHM^\bT_X(Z))$ of a equivariant mixed Hodge module (complex) on $Z$ with support in $X=X_P$, 
induces a weight function on the faces of $P$ defined by $f_Q(y)=\chi_y(i_{x_Q}^*[\cM])$, with $x_Q\in O_{\sigma_Q}$ a point in the orbit of $X_P$ associated to (the cone $\sigma_Q$ in the inner normal fan of $P$ corresponding to) the face $Q$ and $i_{x_Q}: \{x_Q\} \to Z$ the inclusion map. Moreover, any weight function $f$ can be obtained in this way from the class of such a mixed Hodge module complex $\cM$, e.g., by using direct sums and shifts of classes of $(j_Q)_! \bQ^H_{O_{\sig_Q}}(n_Q)$, with $n_Q \in \bZ$ and $j_Q:O_{\sig_Q} \hookrightarrow Z$ the inclusion, regarded as $\bT$-equivariant objects on $Z$  (or, using the subgroup $\chi^\bT_{Hdg}\left(\tau(F^{\bT}(X)[y^{\pm 1}])\right)$ as in the Introduction). In particular, the expressions $E^\varphi_{P,f}(\ell, y)$ and $\wti{E}^\varphi_{P,f}(\ell, y)$ are then also obtained by evaluating a {\it polynomial} $E^\varphi_{P,f}(z, y)$, resp., $\wti{E}^\varphi_{P,f}(z, y)$ at $z=\ell$. \qed
\er

Before studying the generalized weighted Ehrhart polynomial 
$E^\varphi_{P,f}(\ell, y)$, resp., $\wti{E}^\varphi_{P,f}(\ell, y)$,  in more detail, let us also mention the following result, similar  
to Theorem \ref{wlpca}.
\bt\label{thmmin} In the notations of Theorem \ref{wlpca}, we have for $\ell \in \bZ_{>0}$:
\be\label{wcab}
\begin{split}
\chi^\bT_y(X,-\ell D;[\cM]) &= 
 \sum_{Q \preceq P} \chi_y(i_{x_Q}^*[\cM]) \cdot (-1-y)^{\dim(Q)} \cdot \sum_{m \in {\ell} Q \cap M} \chi^{m}.
\end{split}
\ee
\et

\begin{proof}
The proof is very similar to that of Theorem \ref{wlpca}, so we only sketch it below. In view of \eqref{drt2g}, this amounts to showing that for any face $Q$ of $P$, and using notations from loc.cit., we have the following formula:
\be\label{lem2a}
\begin{split}
\chi^\bT(X,(k_Q)_*[\omega_{X_Q}]_\bT \otimes \cO_{X}(-\ell D))&=a_* \left( (k_Q)_*[\omega_{X_Q}]_\bT \otimes \cO_{X}(-\ell D) \right) \\ &= (-1)^{\dim(Q)} \cdot \sum_{m \in \ell Q \cap M} \chi^{m}.
\end{split}
\ee
Indeed,
\begin{equation*}
\begin{split}
a_* \left( (k_Q)_*[\omega_{X_Q}]_\bT \otimes \cO_{X}(-\ell D) \right) 
 &= a_*  (k_Q)_* \left([\omega_{X_Q}]_\bT \otimes k_Q^*\cO_{X}(-\ell D) \right) \\
&\overset{(\ref{tran})}{=} a_*  \left([\omega_{X_Q}]_\bT \otimes \cO_{X_Q}(-\ell D_{Q_0}+div(\chi^{\ell m_0})) \right) \\
&=\chi^{\ell m_0} \cdot a_*\left( [\omega_{X_Q}]_\bT \otimes \cO_{X_Q}(-\ell D_{Q_0}) \right) \\
&\overset{({\ast})}{=} 
\chi^{\ell m_0} \cdot a_*\left( [\omega_{X_Q}]_{\bT_Q} \otimes \cO_{X_Q}(-\ell D_{Q_0}) \right) \\
&\overset{(\ref{am56})}{=} 
\chi^{\ell m_0} \cdot \left(   (-1)^{\dim(Q)} \cdot 
\sum_{m \in \ell (Q-m_0) \cap M} \chi^{m}  \right)\\
&=(-1)^{\dim(Q)} \cdot \sum_{m \in \ell Q \cap M} \chi^{m},
 \end{split}
\end{equation*}
where $(\ast)$  uses again the identification \eqref{equality}. 
\end{proof}

We can now justify the reciprocity formula \eqref{rei} without relying on Brion--Vergne's combinatorial work \cite{BV}, as well as the constant term formula \eqref{ctt}:

\bc Let $P\subset M_\bR \cong \bR^n $ be a full-dimensional lattice polytope, with $\varphi:M_\bR \cong \bR^n \to \bC$ a homogeneous polynomial function. 
For any Laurent polynomial weight vector $\{f_Q\}_{Q \preceq P}$ defined on the set of faces of $P$,  
we have
\be\label{recip} E^\varphi_{P,f}(-\ell, y)=\sum_{Q \preceq P} f_Q(y) \cdot (-1-y)^{\dim(Q)+\deg(\varphi)} \cdot \sum_{m \in  {{\ell}} Q \cap M} \varphi(m)\ee
and 
\be\label{recipt} \wti{E}^\varphi_{P,f}(-\ell, y)=\sum_{Q \preceq P} f_Q(y) \cdot (-1-y)^{\dim(Q)} \cdot (-1)^{\deg(\varphi)} \cdot \sum_{m \in  {{\ell}} Q \cap M} \varphi(m).\ee
In particular,
\be\label{ctt2}E^\varphi_{P,f}(0, y)=\sum_{Q \preceq P} f_Q(y) \cdot (-1-y)^{\dim(Q)+\deg(\varphi)} \cdot \varphi(0)\ee
and 
\be\label{ctt2t}\wti{E}^\varphi_{P,f}(0, y)=\sum_{Q \preceq P} f_Q(y) \cdot (-1-y)^{\dim(Q)} \cdot (-1)^{\deg(\varphi)} \cdot \varphi(0).\ee
\ec

\begin{proof}
As in Remark \ref{r34}, choose a $[\cM]\in K_0(\MHM^\bT_X(Z))$ with $f_Q(y)=\chi_y(i_{x_Q}^*[\cM])$, such that $ E^\varphi_{P,f}(\ell, y)= E^\varphi_{P,\cM}(\ell, y)$. From to the definition of $ E^\varphi_{P,\cM}(\ell, y)$  in \eqref{whpd}, $E^\varphi_{P,\cM}(-\ell, y)$ is obtained from $\chi^\bT_y(X,-\ell D;\cM)$ by substituting each $\chi^{m}$ by 
$$\chi^{m} \mapsto \varphi\left((1+y)\cdot (-m)\right)=(-1-y)^{\deg(\varphi)}\cdot \varphi(m).$$ Formula \eqref{recip} follows then by applying this rule to the expression for $\chi^\bT_y(X,-\ell D;\cM)$ given by \eqref{wcab}. The constant term formula \eqref{ctt2} is then an easy consequence, by setting $\ell=0$ in \eqref{recip}, or by combining \eqref{zero} and \eqref{degg}. Here, we are of course using the fact that $E^\varphi_{P,f}(-\ell, y)$ is a polynomial in $\ell$, a fact proved later on in Theorem \ref{poli}, but see also \cite[Prop.4.1]{BV}. Similarly, for $\wti{E}^\varphi_{P,f}(-\ell, y)$ and $\wti{E}^\varphi_{P,f}(0, y)$. 
\end{proof}

\bex\label{ex65}
Let us next give some examples for Theorem \ref{wlpca}, formulated in terms of our generalized weighted Ehrhart polynomials $E^\varphi_{P,\cM}(\ell, y)$ via Definition \ref{def65}. Let $X=X_P$ be  the projective toric variety with inner normal fan $\Sig=\Sig_P$ and ample Cartier divisor $D=D_P$ associated to a full-dimensional lattice polytope $P\subset M_{\bR} \cong \bR^n$, and let $\ell \in \bZ_{>0}$. For a face $Q$ of $P$ denote by $j_Q:O_{\sig_Q} \hookrightarrow Z$ the inclusion of the orbit  associated to the (cone of the) face $Q$.
\begin{enumerate}
\item[(a)] For $\cM=(j_Q)_! \bQ^H_{O_{\sig_Q}}$, we get
\be\label{69} E^\varphi_{P,\cM}(\ell, y)=(1+y)^{\dim(Q)+\deg(\varphi)} \cdot \sum_{m \in  \Relint({{\ell}} Q) \cap M} \varphi(m),\ee
with $E^\varphi_{P,\cM}(0, y)=(-y-1)^{\dim(Q)+\deg(\varphi)} \cdot \varphi(0)=(-y-1)^{\deg(\varphi)}\cdot \chi_y(O_{\sig_Q}) \cdot \varphi(0)$.

In particular, by evaluating $E^\varphi_{P,\cM}(\ell, y)$ at $y=0$, we get the polynomial
$$E^\varphi_{P,\cM}(\ell, 0)= \sum_{m \in  \Relint({{\ell}} Q) \cap M} \varphi(m)$$ summing up the values of $\varphi$ at lattice points in the relative interior of the dilated face $\ell Q$ (i.e., the Euler-Maclaurin sum for the relative interior of $\ell Q$), with 
$E^\varphi_{P,\cM}(0, 0)=(-1)^{\dim(Q)+\deg(\varphi)}\cdot \varphi(0)$.

\item[(b)] Let $X':=X_{P'}$ be a torus-invariant closed algebraic subset of $X=X_P$ corresponding to a polytopal subcomplex $P' \subseteq P$ (i.e., a closed union of faces of $P$), with $[\cM]=[\bQ_{X'}^H]$ pushed forward to the ambient $\bT$-manifod $Z$ containing $X_P$ as a $\bT$-invariant closed subvariety. Then 
\be
E^\varphi_{P,\cM}(\ell, y)= \sum_{Q \preceq P'} (1+y)^{\dim(Q)+\deg(\varphi)} \cdot \sum_{m \in  \Relint({{\ell}} Q) \cap M} \varphi(m)  \in \bZ[y],
\ee
where the summation is over the faces $Q$ of $P'$. Hence, 
 $$E^\varphi_{P,\cM}(\ell, 0)=\sum_{m\in \ell P' \cap M} \varphi(m)$$ is the corresponding Euler--Maclaurin sum for $\ell P'$, and 
$$E^\varphi_{P,\cM}(0, y)=\sum_{Q \preceq P'} (-y-1)^{\dim(Q)+\deg(\varphi)}\cdot \varphi(0)=(-y-1)^{\deg(\varphi)} \cdot \chi_y(X')\cdot \varphi(0),$$
which gives that $$E^\varphi_{P,\cM}(0, 0)=(-1)^{\deg(\varphi)} \cdot  \chi(P') \cdot \varphi(0),$$ with $\chi(P')$  the topological Euler characteristic of $P'$. 

\item[(c)] For $\cM=(j_Q)_* \bQ^H_{O_{\sig_Q}}$, we get by \eqref{drtz} and  \eqref{lem2b} that 
\be\label{pu} E^\varphi_{P,\cM}(\ell, y)=(1+y)^{\dim(Q)+\deg(\varphi)} \cdot \sum_{m\in \ell Q \cap M} \varphi(m) ,\ee
with $E^\varphi_{P,\cM}(0, y)=(1+y)^{\deg(\varphi)} \cdot   \chi_y(X; [(j_Q)_* \bQ^H_{O_{\sig_Q}}]) \cdot \varphi(0)$, $$E^\varphi_{P,\cM}(\ell, 0)=\sum_{m\in \ell Q \cap M} \varphi(m)$$ the Euler--Maclaurin sum for $\ell Q$, and $E^\varphi_{P,\cM}(0, 0)=\chi(Q)\cdot \varphi(0)$.
\qed
\end{enumerate}
\eex

Let $X=X_P$ be as before the projective toric variety with inner normal fan $\Sig=\Sig_P$ and ample Cartier divisor $D=D_P$ associated to a full-dimensional lattice polytope $P\subset M_{\bR} \cong \bR^n$. Let $X \hookrightarrow Z$ be the equivariant embedding of $X$ into an ambient $\bT$-manifold, as in \cite{Su}.
Then another important example is provided by $\cM={IC'}^H_{V_{Q'}}$ with $V_{Q'}$ the closure of the orbit corresponding to the face $Q'$ of $P$.  This is a shifted $\bT$-equivariant pure Hodge module on $Z$ with support in $V_{Q'} \subset X$.
As explained in Example \ref{g-face} of the Introduction, 
the stalk contributions $\chi_y({IC'}_{V_{Q'}}^H\vert_{x_Q})$ have the following combinatorial description (see, e.g., \cite{BM,F, DL,Sa}, etc):
\be\label{co}
\chi_y({IC'}_{V_{Q'}}^H\vert_{x_Q})=\wti{g}_Q(-y):=g_{Q^\circ}(-y),
\ee
with $Q$ a face of $Q'$ and $Q^\circ$  the corresponding face of a polar polytope ${P'}^\circ$ as explained in Example \ref{g-face}, and $g_{Q^\circ}(t) \in \bZ[t]$ the Stanley {\it $g$-polynomial} of $Q^\circ$, cf. \cite[Sect.3]{St}.
In the notations of Definition \ref{def65}, we then get for  $E^{\varphi}_{Q',IC'}(\ell, y):= E^{\varphi}_{P,{IC'}^H_{V_{Q'}}}(\ell, y)$ the following.

\bc\label{wlpc}
Let $X=X_P$ is the projective toric variety with ample Cartier divisor $D=D_P$ associated to a full-dimensional lattice polytope $P\subset M_{\bR} \cong \bR^n$, with $\varphi:M_\bR \cong \bR^n \to \bC$ a homogeneous polynomial function. Let $Q'$ be a fixed face of $P$. Then, in the above notations, one has for $\ell \in \bZ_{>0}$, 
\be\label{wc}
\begin{split}
E^{\varphi}_{Q',IC'}(\ell, y)&= \sum_{Q \preceq Q'} \wti{g}_Q(-y) \cdot (1+y)^{\dim(Q)+\deg(\varphi)} \cdot \sum_{m \in  \Relint({{\ell}} Q) \cap M} \varphi(m), 
\end{split}
\ee
and
\be\label{recip2} E^\varphi_{Q',IC'}(-\ell, y)=\sum_{Q \preceq Q'} \wti{g}_Q(-y) \cdot (-1-y)^{\dim(Q)+\deg(\varphi)} \cdot \sum_{m \in  {{\ell}} Q \cap M} \varphi(m).\ee
where the summation is over the faces $Q$ of $P$.
In particular,
\be\label{ctt22}E^\varphi_{Q',IC'}(0, y)=\sum_{Q \preceq Q'} \wti{g}_Q(-y) \cdot (-1-y)^{\dim(Q)+\deg(\varphi)} \cdot \varphi(0).\ee
 Similar formulas hold for $\wti{E}^{\varphi}_{Q',IC'}(\ell, y)$ by deleting the factor $(1+y)^{\deg(\varphi)}$.
\ec

\br \begin{itemize}
\item[(a)] In the case when $Q'=P$, formula \eqref{wc} equals $G_{\varphi}(\ell,y)$ from \cite[formula (14)]{BGM}.
\item[(b)]  If the full-dimensional polytope $P$ is simple, then $X=X_P$ is a simplicial toric variety, hence $[{IC'}_{V_{Q'}}^H] = [\bQ^H_{V_{Q'}}] \in K_0^\bT(\MHM_X(Z))$. So, in this case,  $\chi_y({IC'}_{V_{Q'}}^H\vert_{x_Q})=\wti{g}_Q(-y)=1$. \end{itemize} 
\qed
\er

We next discuss the effect of duality on generalized weighted Ehrhart polynomials, proving a generalized {\it reciprocity formula}. This takes the form of a {\it purity} statement for pure equivariant Hodge modules (like $IC_{V_{Q'}}^H$ for $Q'$ a face of $P$).

\bt[Reciprocity and Purity]\label{repr} In the above notations, for any Grothendieck class
$[\cM]\in \im(\chi^\bT_{Hdg})\subset K_0(\MHM^\bT_X(Z))$, we have
\be\label{r1}
\begin{split}
E^\varphi_{P,\cM}(-\ell, y)=(-y)^{\deg(\varphi)} \cdot E^\varphi_{P, \bD^\bT_Z\cM}(\ell, \frac{1}{y}),
\end{split}
\ee
and 
\be\label{r1t}
\wti{E}^\varphi_{P,\cM}(-\ell, y)=(-1)^{\deg(\varphi)} \cdot \wti{E}^\varphi_{P, \bD^\bT_Z\cM}(\ell, \frac{1}{y}).
\ee
In particular, if $\cM$ is such a self-dual pure equivariant Hodge module of weight $n'$ on $X=X_P$,
then the following purity property holds:
\be\label{r1b}
E^\varphi_{P,\cM}(-\ell, y)=(-y)^{n'+\deg(\varphi)} \cdot E^\varphi_{P, \cM}(\ell, \frac{1}{y}),
\ee
and \be\label{r1bt}
\wti{E}^\varphi_{P,\cM}(-\ell, y)=(-y)^{n'} \cdot (-1)^{\deg(\varphi)} \cdot \wti{E}^\varphi_{P, \cM}(\ell, \frac{1}{y}).
\ee
\et

Before proving Theorem \ref{repr}, we need the following.
\bl\label{laj}
Under the assumptions of Theorem \ref{repr} and for $\ell \in \bZ_{>0}$, we have:
\be\label{chid}
\chi^\bT_y(X,\ell D;\bD^\bT_Z[\cM])=\sum_{Q \preceq P} f_Q(\frac{1}{y}) \cdot (-y)^{-\dim(Q)} \cdot (1+y)^{\dim(Q)} \cdot  \sum_{m \in  {{\ell}} Q \cap M} \chi^{-m},
\ee
with $f_Q(y)=\chi_y(i^*_{x_Q}[\cM])$.
\el

\begin{proof}
First, by \eqref{drt2g} and \eqref{drtc}, we have
\be\label{e3}
\begin{split}
\DR^\bT_y([\cM])&=\sum_{Q \preceq P} \, \chi_y(i^*_{x_Q}[\cM]) \cdot (1+y)^{\dim(Q)} \cdot (k_Q)_* [\omega_{V_{\sigma_Q}}]{_\bT} \\
&=\sum_{Q \preceq P} \, f_Q(y) \cdot  \DR^\bT_y \left([(j_Q)_! \bQ^H_{O_{\sig_Q}}] \right),
\end{split}
\ee
with $j_Q:O_{\sig_Q} \hookrightarrow Z$ and $k_Q:V_{\sig_Q}=X_Q \hookrightarrow Z$ the inclusion maps. Moreover, as in \eqref{e1}, we have:
\be\label{e4}\DR^\bT_y \left(\bD^\bT_Z \left([(j_Q)_! \bQ^H_{O_{\sig_Q}}]\right) \right)
= (-y)^{-\dim(Q)} \cdot   \DR^\bT_y \left([(j_Q)_* \bQ^H_{O_{\sig_Q}}]\right). \ee
Since $\DR^\bT_y$ commutes with dualities, up to exchanging $y$ by $1/y$, we then get from \eqref{e3} and \eqref{e4} that
\be\label{e5}
\begin{split}
\DR^\bT_{y}(\bD^\bT_Z[\cM])&=\sum_{Q \preceq P} \, f_Q(\frac{1}{y}) \cdot  (-y)^{-\dim(Q)} \cdot \DR^\bT_{y}\left([(j_Q)_* \bQ^H_{O_{\sig_Q}}]\right)\\
&\overset{\eqref{drtz}}{=}\sum_{Q \preceq P} \, f_Q(\frac{1}{y}) \cdot  (-y)^{-\dim(Q)} \cdot (1+y)^{\dim(Q)} \cdot (k_Q)_*[\cO_{X_Q}]_\bT.
\end{split}\ee

Therefore, with $a$ denoting as before the constant map to a point,
\be
\begin{split}
\chi^\bT_y(X,\ell D;\bD^\bT_Z[\cM])&=a_*\left(\DR^\bT_y(\bD^\bT_Z[\cM])  \otimes \cO_X(\ell D)  \right) \\
&=\sum_{Q \preceq P} \, f_Q(\frac{1}{y}) \cdot  (-y)^{-\dim(Q)} \cdot (1+y)^{\dim(Q)} \cdot a_* \left( (k_Q)_*[\cO_{X_Q}]_\bT \otimes \cO_{X}(\ell D) \right) \\
&\overset{\eqref{lem2b}}{=}
\sum_{Q \preceq P} \, f_Q(\frac{1}{y}) \cdot  (-y)^{-\dim(Q)} \cdot (1+y)^{\dim(Q)} \cdot \sum_{m \in \ell Q \cap M} \chi^{-m}.
\end{split}\ee
\end{proof}

\begin{proof}[Proof of Theorem \ref{repr}]
First, we get from \eqref{chid} that
$$
E^\varphi_{P, \bD^\bT_Z\cM}(\ell, {y})=\sum_{Q \preceq P} \, f_Q(\frac{1}{y}) \cdot  (-y)^{-\dim(Q)} \cdot (1+y)^{\dim(Q)} \cdot (1+y)^{\deg(\varphi)} \cdot \sum_{m \in \ell Q \cap M} \varphi(m),
$$
with $f_Q(y)=\chi_y(i^*_{x_Q}[\cM])$ as before. Substituting $1/y$ in place of $y$, we then get
\be
\begin{split}
E^\varphi_{P, \bD^\bT_Z\cM}(\ell, \frac{1}{y})&=\sum_{Q \preceq P} \, f_Q({y}) \cdot  (-y)^{\dim(Q)} \cdot (1+\frac{1}{y})^{\dim(Q)+\deg(\varphi)}  \cdot \sum_{m \in \ell Q \cap M} \varphi(m)\\
&=\sum_{Q \preceq P} \, f_Q({y}) \cdot  y^{-\deg(\varphi)} \cdot (-1)^{\dim(Q)} \cdot (1+{y})^{\dim(Q)+\deg(\varphi)}  \cdot \sum_{m \in \ell Q \cap M} \varphi(m)\\
&=\sum_{Q \preceq P} \, f_Q({y}) \cdot  y^{-\deg(\varphi)} \cdot (-1)^{\deg(\varphi)} \cdot (-1-{y})^{\dim(Q)+\deg(\varphi)}  \cdot \sum_{m \in \ell Q \cap M} \varphi(m)\\
&=(-y)^{-\deg(\varphi)} \cdot \sum_{Q \preceq P} f_Q(y) \cdot (-1-y)^{\dim(Q)+\deg(\varphi)} \cdot \sum_{m \in  {{\ell}} Q \cap M} \varphi(m)\\
&\overset{\eqref{recip}}{=} (-y)^{-\deg(\varphi)} \cdot  E^\varphi_{P,f}(-\ell, y).
\end{split}
\ee

Formula \eqref{r1b} follows immediately by combining \eqref{r1} with \eqref{ttw}.

The formulas for the polynomials $\wti{E}$ are obtained similarly.
\end{proof}

\bex
The reciprocity formula \eqref{r1} can be applied to all cases of Example \ref{ex65}. In the (dual) cases $(a)$ and $(c)$, for $\cM=(j_Q)_! \bQ^H_{O_{\sig_Q}}$ with $\bD^\bT_Z [\cM]=[(j_Q)_* \bQ^H_{O_{\sig_Q}}(\dim(Q))]\in K^\bT_0(\MHM_X(Z))$, formula \eqref{r1} specializes to the {\it classical Ehrhart reciprocity}  for $Q$ from \cite[Prop.4.1]{BV}, compare with the examples after Theorem \ref{thm-pur} in the Introduction. Indeed, 
$$
E^\varphi_{P,\cM}(\ell, y)=(1+y)^{\dim(Q)+\deg(\varphi)} \cdot \sum_{m\in \Relint({{\ell}} Q) \cap M} \varphi(m)
$$
and 
$$
 (-y)^{\deg(\varphi)} \cdot E^\varphi_{P,\bD^\bT_X\cM}(\ell, \frac{1}{y})=(1+y)^{\dim(Q)+\deg(\varphi)} \cdot \left(  (-1)^{\dim(Q)+\deg(\varphi)} \sum_{m\in {{\ell}} Q \cap M} \varphi(m) \right).
$$ \qed
\eex

\bex[$IC$-reciprocity]\label{icre}
Let $Q'$ be a face of $P$ with $V_{Q'}$ the closure  of the corresponding torus orbit in $X$. Let $\cM={IC'}^H_{V_{Q'}}$.  
As explained in Example \ref{g-face} of the Introduction, the corresponding  weight function for a face $Q$ of $Q'$ is given by Stanley's $g$-function  $f_Q(y)=\wti{g}_Q(-y)$ of the corresponding face of a polar polytope for $Q'$. Then our general purity formula \eqref{r1b} recovers via Corollary \ref{wlpc} the reciprocity formula proved  in \cite[Theorems 1.3 and 2.6]{BGM} by combinatorial methods for $Q'=P$.\qed
\eex

 \appendix
 
\section{}
In this final section we give a geometric proof of the following result (without relying on the combinatorial work of Brion--Vergne \cite{BV}).

\bt\label{poli}
Let $P\subset M_R \cong \bR^n$ be a full-dimensional lattice polytope. Let $\varphi: M_\bR \cong \bR^n \to \bC$ be a homogeneous polynomial function, and $\ell \in \bZ_{>0}$. Then for any Laurent polynomial {weight function} $$f:\{\text{faces of } P\} \to \bZ[y^{\pm 1}], \ \ \ Q \mapsto f_Q(y)$$ defined on the set of faces of $P$, the expression 
\be\label{whpga} E^\varphi_{P,f}(\ell, y):=\sum_{Q \preceq P} f_Q(y) \cdot (1+y)^{\dim(Q)+\deg(\varphi)} \cdot \sum_{m \in  \Relint({{\ell}} Q) \cap M} \varphi(m)= (1+y)^{\deg(\varphi)} \cdot \wti{E}^\varphi_{P,f}(\ell, y) \ee
is obtained by evaluating a {polynomial} $E^\varphi_{P,f}(z, y)$ at $z=\ell$.
\et

Let $X=X_P$ be the toric variety with torus $\bT$ and ample Cartier divisor $D=D_P$ associated to the lattice polytope $P$. Let $Z$ be an ambient smooth $\bT$-variety.
In view of Remark \ref{r34}, it suffices to prove the assertion in the above theorem for $E^\varphi_{P,\cM}(\ell, y)$, as defined in \eqref{whpd}, i.e., assuming that $f_Q(y)=\chi_y(i_{x_Q}^*[\cM])$, for a certain class $[\cM]\in K_0(\MHM^\bT_X(Z))$  of a mixed Hodge module (complex) on $Z$ with support in $X=X_P$,  with $x_Q\in O_{\sigma_Q}$ a point in the orbit of $X$ associated to (the cone $\sigma_Q$ in the inner normal fan of $P$ corresponding to) the face $Q$, and $i_{x_Q}: \{x_Q\} \to Z$ the inclusion map.
Our proof in this case is based on localization techniques developed in \cite{CMSSEM}. We first introduce some terminology.

Recall that the (Borel-type) rational equivariant cohomology of $X$ is defined as \be H^*_\bT(X;\bb{Q}):=H^*(E\bT\times_\bT X;\bb{Q}),\ee
where $\bT\hookrightarrow E\bT\to B\bT$ is the universal principal $\bT$-bundle, i.e., $B\bT=(B\bb{C}^*)^n=(\bC P^{\infty})^n$ and $E\bT$ is contractible. In particular, if $X=pt$ is a point space, 
\be H^*_\bT(pt;\bQ)=H^*(B\bT;\bQ)\simeq \bQ[t_1,\ldots,t_n]=:(\Lambda_\bT)_{\bQ}.\ee
Similarly, one can define equivariant homology groups $H_*^\bT(X;\bb{Q})$ as in \cite{BZ},  \cite[Sect.17]{AF} and \cite[Sect.2.5]{CMSSEM}, using suitable finite-dimensional approximations of $B\bT$. Note that these equivariant homology groups can be non-zero also in negative degrees. As in the non-equivariant context, there is a cap product operation
$$\cap: H^*_\bT(X;\bb{Q}) \times H_*^\bT(X;\bb{Q})  \lra H_*^\bT(X;\bb{Q}),$$
together with an equivariant Poincar\'e duality if $X$ is smooth, see \cite[Sect.17.1]{AF}. In particular, $H^*_\bT(X;\bb{Q})$, resp., $H_*^\bT(X;\bb{Q})$, is a $(\Lambda_\bT)_{\bQ}$-algebra, resp., $(\Lambda_\bT)_{\bQ}$-module.

For $m\in M$, an element in the character lattice of the torus $\bT$, one can view the corresponding character $\chi^m:\bT \to \bC^*$ as a $\bT$-equivariant line bundle $\bC_{\chi^m}$ over a point space $pt$, where the $\bT$-action on $\bC$ is induced via $\chi^m$.
This gives rise to an abelian group isomorphism
$$(M,+) \simeq (Pic_\bT(pt), \otimes) \:.
$$
Recall that $m\mapsto -m$ corresponds to the duality involution $(-)^{\vee}$. Taking the first equivariant Chern class 
$c^1_\bT$ (or the dual $-c^1_\bT=c^1_\bT\circ (-)^{\vee}$) of $\bC_{\chi^m}$ gives an  isomorphism
\begin{equation}\label{1}
c=c^1_\bT,  \:\text{resp.}, \:  s=-c^1_\bT :\:
M \simeq H_\bT^2(pt;\bb{Z})
\end{equation}
and
\begin{equation}\label{2}
 c,  \:\text{resp.}, \: s: Sym_{\bb{Q}}(M)  \simeq  H_\bT^*(pt;\bb{Q})=(\Lambda_\bT)_{\bb{Q}}\:, 
\end{equation}
with $Sym_{\bb{Q}}(M)=\bigoplus_{k=0}^{\infty} Sym^k_{\bb{Q}}(M)$ the (rational) symmetric algebra of $M$.
So if $m_i$ ($i=1,\dots,n$) is a basis of $M\simeq \bb{Z}^n$, then  $H_\bT^*(pt;\bb{Q})=(\Lambda_\bT)_{\bb{Q}}\simeq \bb{Q}[t_1,\dots,t_n],$  with
$t_i=\pm c^1_\bT(\bC_{\chi^{m_i}})$ for $i=1,\dots,n$. 

Let $$ \bb{Q}\{t_1,\dots,t_n\}\simeq (H^*_\bT(pt;\bb{Q}))^{an}=:(\Lambda^{an}_\bT)_{\bb{Q}}\subset (\widehat{\Lambda}_\bT)_{\bb{Q}}\simeq \bQ[[t_1,\ldots, t_n]]$$ be the subring of {convergent power series} (around zero) with rational coefficients, i.e., after pairing with $z\in N_{\bb{K}}=N\otimes_{\bb{Z}}\bb{K}$ (for $\bb{K}=\bb{R},\bb{C}$) one gets a convergent power series {\it function} in $z$ around zero, whose corresponding Taylor polynomials have rational coefficients. Here, the degree completion $(\widehat{\Lambda}_\bT)_{\bb{Q}}$ of $(\Lambda_\bT)_{\bb{Q}}$ corresponds to formal power series with rational coefficients.

Recall that if $\cF$ is a  $\bT$-equivariant coherent sheaf on $X$, then the  {$K$-theoretic Euler characteristic} of $\cF$ is  defined as in \eqref{Kchi} by: \be\label{Kchia}\begin{split} \chi^\bT(X,\cF)&=\sum_{m \in M} \sum_{i=0}^{\dim(X)} (-1)^i \dim_\bC H^i(X;\cF)_{\chi^m} \cdot \chi^m \in \bZ[M], \end{split} \ee
By applying the equivariant Chern character ring homomorphism $$\ch^\bT:\bZ[M] \simeq K^0_\bT(pt) \hookrightarrow  (\Lambda^{an}_\bT)_{\bb{Q}}\subset 
 (\widehat{\Lambda}_\bT)_{\bb{Q}}$$
 given by $m \mapsto e^{c(m)}$, 
  we further get the  {\it cohomological Euler characteristic} of $\cF$:
\be\label{f41}
\chi^\bT(X,\cF)=\sum_{m \in M} \sum_{i=0}^{\dim(X)} (-1)^i \dim_\bC H^i(X;\cF)_{\chi^m} \cdot e^{c(m)} \in (\Lambda^{an}_\bT)_{\bb{Q}}
\subset  (\widehat{\Lambda}_\bT)_{\bb{Q}}\:.
\ee

Consider now the fixed-point set $X^\bT$ of the torus action; this can be identified with the set of vertices of the lattice polytope $P$. 
By \cite[Lem.8.4, Lem.8.5]{BZ}  and \cite[Thm.17.3.1]{AF}, the inclusion $\alpha: X^\bT \hookrightarrow X$ induces an injective homomorphism of $H^{*}_\bT(pt;\bQ)$-modules $$\alpha_*:\widehat{H}_{*}^\bT(X^\bT;\bQ)  \hookrightarrow \widehat{H}_{*}^\bT(X;\bQ)$$
which becomes an isomorphism $$\alpha_*:\widehat{H}_{*}^\bT(X^\bT;\bQ)_L  \simeq \widehat{H}_{*}^\bT(X;\bQ)_L$$
upon localization at the multiplicative subset $L\subset (\Lambda_\bT)_\bQ=H^{*}_\bT(pt;\bQ)$ generated by the elements $\pm c(m)$, for $0\neq m \in M$ (cf. \cite[Cor.8.9]{BZ}). Here, we use the completion  of 
$$ {H}_{*}^\bT(X; \bQ) :=\bigoplus_{i \leq 2\dim(X)} {H}_{i}^\bT(X;\bQ),$$
defined as
$$ \widehat{H}_{*}^\bT(X; \bQ) :=\prod_{i \leq 2\dim(X)} {H}_{i}^\bT(X;\bQ).$$ 
 Since $$\widehat{H}_{*}^\bT(X^\bT;\bQ) = \bigoplus_{x \in X^\bT} \widehat{H}_{*}^\bT(x;\bQ),$$
one gets via the isomorphism $\alpha_*$ a projection map of $H^{*}_\bT(pt;\bQ)$-modules, called  the {homological localization map at $x$}, 
$$pr_x:\widehat{H}_{*}^\bT(X;\bQ)_L \simeq  \bigoplus_{x \in X^\bT} \widehat{H}_{*}^\bT(x;\bQ)_L \lra \widehat{H}_{*}^\bT(x;\bQ)_L = L^{-1}(\widehat{\Lambda}_\bT)_\bQ.$$
Moreover, as in \cite[Sect.4.2, formula (4.22)]{CMSSEM}, if $i_x:\{x\} \hookrightarrow X$ is the inclusion map, then for any $a \in \widehat{H}^{*}_\bT(X;\bQ)_L$ and $b \in \widehat{H}_{*}^\bT(X;\bQ)_L$, we have:
\be\label{proj}
pr_x(a \cap b)=i_x^*(a) \cap pr_x(b).
\ee

We then have the following result from \cite[Sect.4.2, Prop.4.19]{CMSSEM}, adapted to the notations of our setting.
\bp
Let $X=X_P$ be the projective toric variety associated to the full-dimensional lattice polytope $P\subset M_\bR\cong \bR^n$, with ample Cartier divisor $D=D_P$. Then for $\cF$ a $\bT$-coherent sheaf on $X$, the cohomological Euler characteristic $\chi^\bT(X,\cO_X(D) \otimes \cF)$ can be calculated via localization at the $\bT$-fixed points as:
\be\label{eqec}
\chi^\bT(X,\cO_X(D) \otimes \cF)=\sum_{v \ {\rm vertex \ of }  \ P} 
 pr_{x_v} (\td_*^\bT([\cO_X(D) \otimes \cF]) ,
 \ee
 with each summand 
\be\label{ana} pr_{x_v} (\td_*^\bT([\cO_X(D) \otimes \cF]) \in   L^{-1}(\Lambda^{an}_\bT)_\bQ ,\ee
where $x_v$ is the $\bT$-fixed point of $X$ corresponding to the vertex $v$ of $P$, and $\td_*^\bT:K^\bT_0(X) \lra \widehat{H}^\bT_{2*}(X;\bQ)$ is the equivariant Riemann--Roch map of Edidin--Graham \cite{EG} (see also \cite{BZ}).
\ep

\br
The above homological localization formula is in fact a consequence of a similar localization formula in equivariant $K$-theory, using the fact that the equivariant Chern character induces a ring homomorphism of localized rings
$$\ch^\bT: S^{-1}\bZ[M]\simeq S^{-1}K^0_\bT(pt) \lra L^{-1}(\Lambda^{an}_\bT)_\bQ,$$
with $S$ the multiplicative system generated by elements $1-\chi^m$ for $0 \neq m \in M$. Note that for a homogeneous polynomial function $\varphi:M_\bR \to \bC$, this property does not hold for the corresponding group homomorphisms $\wti{E}^\varphi$ and $E^\varphi$ as defined in the Introduction. For this reason, we use exponential sums and their Taylor expansions (as in \cite{BV}), as we shall explain below.
\qed
\er

Recall the identification
$s: Sym_{\bb{Q}}(M)  \simeq  H_\bT^*(pt;\bb{Q})=:(\Lambda_\bT)_{\bb{Q}}$. 
Let us fix $z \in N_\bC:=N \otimes_\bZ \bC= \Hom_\bR(M_\bR, \bC)$. By the universal property of $Sym$, $z$ induces  an $\bR$-algebra homomorphism
$$\langle -, z \rangle : Sym_{\bb{R}}(M) \to \bC \:,$$
by which we can view $\langle p, z \rangle$ for $z$ now variable and $p\in  Sym_{\bb{R}}(M)$, resp., $p\in  (Sym_{\bb{R}}(M))^{an}$ fixed,
as a $\bC$-valued polynomial on $N_\bR$, resp., as a {convergent  power series function} around zero in  $N_\bR$.

We can now complete the proof of Theorem \ref{poli}.

\begin{proof}[Proof of Theorem \ref{poli}] 

From the geometric viewpoint, it is more natural to prove the assertion for  $\wti{E}^\varphi_{P,f}(\ell, y)$.

By applying the Chern character $\ch^\bT$ to formula \eqref{wca}, we obtain its cohomological counterpart as follows:
\be\label{chic}
\begin{split}
\sum_{Q \preceq P} f_Q(y) \cdot (1+y)^{\dim(Q)} \cdot \sum_{m \in \Relint({{\ell}} Q) \cap M} e^{s(m)}&=
\chi^\bT_y(X,\ell D;[\cM]) \\
&=  
\chi^\bT(X;\DR^\bT_y([\cM]) \otimes \cO_X(\ell D))\\
&\overset{\eqref{eqec}}{=} \sum_{v \ {\rm vertex}} 
 pr_{x_v} (\td_*^\bT(\DR^\bT_y([\cM]) \otimes \cO_X(\ell D)) ) \\
 &\overset{(*)}{=} \sum_{v \ {\rm vertex}} 
 pr_{x_v}  (e^{\ell [D]} \cap \td_*^\bT(\DR^\bT_y([\cM]) )) \\
 &\overset{\eqref{proj}}{=} \sum_{v \ {\rm vertex}} 
e^{\ell \cdot i_{{x_v}}^*[D]} \cap  pr_{x_v}  ( \td_*^\bT(\DR^\bT_y([\cM]) )),
\end{split}
\ee
where $i_{{x_v}}: \{x_v\} \hookrightarrow X$ is the inclusion of the $\bT$-fixed point of $X$ corresponding to the vertex $v$ of $P$,  
with $f_Q(y)=\chi_y(i_{x_Q}^*[\cM])$, and $(*)$ makes use of the module property for $\td_*^\bT$.

Applying $\langle -, z \rangle$ to both sides of \eqref{chic}, we obtain 
\be\label{chic2}
\begin{split}
\sum_{Q \preceq P} f_Q(y) \cdot (1+y)^{\dim(Q)} \cdot \sum_{m \in \Relint({{\ell}} Q) \cap M} e^{\langle m,z \rangle}&=
\sum_{v \ {\rm vertex}} 
e^{\langle \ell \cdot i_{{x_v}}^*[D] , z \rangle} \cdot  \langle pr_{x_v}  ( \td_*^\bT(\DR^\bT_y([\cM]) )), z \rangle.
\end{split}
\ee

Let us now fix a vertex $v$ of $P$. By \eqref{ana}, one also has that $$pr_{x_v}  ( \td_*^\bT(\DR^\bT_y([\cM]) )) \in L^{-1}(\Lambda^{an}_\bT)_\bQ[y^{\pm 1}],$$
since $e^{\ell \cdot i_{{x_v}}^*[D]} \in (\Lambda^{an}_\bT)_\bQ$ is invertible. Then $\langle pr_{x_v}  ( \td_*^\bT(\DR^\bT_y([\cM]) )), z \rangle$ is a finite linear combination over $\bZ[y^{\pm 1}]$ of terms of the form $\frac{h(z)}{\prod_{i=1}^k \langle m_i , z \rangle}$, with $0 \neq m_i \in M$ and $h$ analytic in $z$ near $0$. 

As in \cite[Prop.4.1]{BV}, we next replace $z$ by $tz$ in \eqref{chic2}, with $0\neq t \in \bR$ small and $z$ generic (in the sense that $\langle m_i , z \rangle \neq 0$, $i=1,\ldots k$), and consider the expansion into Laurent series in $t$ on both sides of the equality.  
We then obtain exactly as in \cite[Prop.4.1]{BV} the polynomial behavior in $\ell$ of $\wti{E}^\varphi_{P,\cM}(\ell, y)$ for the polynomial function $\varphi(m)=\langle m, z \rangle^r$, for any given non-negative integer $r$ and $z$ outside a finite number of hyperplanes. This then implies the statement for any polynomial $\varphi(m)$.
\end{proof}

\br
If we want to estimate the degree in $\ell$ of the generalized weighted Ehrhart polynomial in Theorem \ref{poli} to be at most $deg(\varphi)+d$ for some $d \in \bZ_{>0}$, we need to control uniformly the number $k \leq d$ of linear factors appearing in the above localization. In fact, we can choose $d=n=\dim(P)$, as we will now explain. By linearity over $\bZ[y^{\pm 1}]$ and using the previous notations, it suffices to work with $\cM=(j_Q)_! \bQ^H_{O_{\sig_Q}}$, for all faces $Q$ of $P$. By using a resolution of singularities of the orbit closure and functoriality, one can further reduce to the case of a smooth projective toric variety $X'$ of dimension at most $n'\leq n$. In this smooth context, a more precise formula for the homological localization $pr_x$ to a torus fixed point $x$ is available via equivariant Poincar\'e duality as (e.g., see \cite[Sect.4.2]{CMSSEM}):
$$pr_x=\frac{i_x^*}{c_{n'}^\bT(N_x)},$$
with $i_x:\{x\} \hookrightarrow X'$ the inclusion map and $c_{n'}^\bT(N_x)$ the top equivariant Chern class of the normal bundle $N_{x}$ of $x$ in $X'$. Taking a character decomposition of $N_x=\bigoplus_{i=1}^{n'} \bC_{\chi^{m_i}}$, we get $c_{n'}^\bT(N_x)=\prod_{i=1}^{n'} c(m_i)$. \qed
\er


\begin{thebibliography}{99}

\bibitem{A} P. Achar, {\it Equivariant mixed Hodge modules}, preprint, Lecture notes from the Clay Math- ematics Institute workshop on Mixed Hodge Modules and Applications (2013), available at https://www.math.lsu.edu/~pramod/docs/emhm.pdf.

\bibitem{A2} P. Achar, {\it Perverse Sheaves and Applications to Representation Theory},
Mathematical Surveys and Monographs, no. 258, American Mathematical Society, Providence, RI, 2021. 

\bibitem{AMSS} P. Aluffi, L.. Mihalcea, J. Sch\"urmann, C. Su, {\it Motivic Chern classes of Schubert cells, Hecke algebras, and applications to Casselman's problem}, Ann. Sci. \'Ec. Norm. Sup\'er. (4) 57 (2024), no. 1, 87--141.

\bibitem{AF} D. Anderson, W. Fulton, {\it Equivariant cohomology in algebraic geometry},
Cambridge Stud. Adv. Math., 210
Cambridge University Press, Cambridge, 2024.

\bibitem{BGM} M. Beck, P. Gunnells, E. Materov, {\it Weighted lattice point sums in lattice polytopes, unifying Dehn--Sommerville and Ehrhart--Macdonald},  Discrete Comput. Geom. 65 (2021), no. 2, 365--384.

\bibitem{Bit}  F. Bittner, {\it The universal Euler characteristic for varieties of characteristic zero}, Compos. Math. 140 (2004), no. 4, 1011--1032.

\bibitem{BM}  T. Braden, R. MacPherson, {\it Intersection homology of toric varieties and a conjecture of Kalai}, Comment. Math. Helv. 74 (1999), no. 3, 442--455.

\bibitem{BSY} J.-P. Brasselet, J. Sch\"urmann, S. Yokura, {\it Hirzebruch classes and motivic Chern classes for singular spaces}, J. Topol. Anal. 2 (2010), no. 1, 1--55.

\bibitem{BZ}  J.L. Brylinski, B. Zhang, {\it Equivariant Todd Classes for toric varieties,} 
arXiv:math/0311318.

\bibitem{BV} M. Brion, M. Vergne, {\it Lattice points in simple polytopes}, J. Amer. Math. Soc. 10 (1997), no. 2, 371--392.

 \bibitem{CMSSEM} S. Cappell, L. Maxim, J. Sch\"urmann, J. Shaneson, {\it  Equivariant toric geometry and Euler-Maclaurin formulae}, Comm. Pure Appl. Math. (published online). DOI:10.1002/cpa.70016.

\bibitem{CG} N. Chriss, V. Ginzburg, {\it Representation theory and complex geometry}, Birkh\"auser Boston, Inc., Boston, MA, 1997.

\bibitem{CLS} D. A. Cox, J. B. Little,  H. K. Schenck, {\it Toric varieties}, Graduate Studies in Mathematics, 124. American Mathematical Society, Providence, RI, 2011.

\bibitem{Dan}  V. I. Danilov,   {\it The geometry of toric varieties}, Russian Math. Surveys 33 (1978), 97--154.

\bibitem{DL} J. Denef, F. Loeser, 
{\it Weights of exponential sums, intersection cohomology, and Newton polyhedra}, 
Invent. Math. 106 (1991), no. 2, 275--294. 

\bibitem{DM} B. Dirks, M. Musta\c{t}\u{a}, {\it The Hilbert series of Hodge ideals of hyperplane arrangements}, J. Singul. 20 (2020), 232--250.

\bibitem{EG} D. Edidin, W. Graham,  {\it Riemann-Roch for equivariant Chow groups},  
Duke Math. J. 102 (2000), no. 3, 567--594.

\bibitem{F} K.-H. Fieseler, {\it Rational intersection cohomology of projective toric varieties}, J. Reine Angew. Math. 413 (1991), 88--98. 

\bibitem{FRW}  L. Feher, R. Rimanyi, A. Weber, {\it Motivic Chern classes and K-theoretic stable envelopes}, Proc. Lond. Math. Soc. (3) 122 (2021), no. 1, 153--189.

\bibitem{Fu2}
W. Fulton, 
{\it Introduction to toric varieties}, 
Annals of Mathematics Studies, 131. Princeton University Press, Princeton, NJ, 1993.

\bibitem{FM} W. Fulton, R. MacPherson, {\it Categorical framework for the study of singular spaces}, Memoirs of the
AMS 243 (1981).

\bibitem{HTT}  R. Hotta, K.  Takeuchi, T. Tanisaki, {\it D-modules, perverse sheaves, and representation theory},
Progr. Math., 236 Birkh\"auser Boston, Inc., Boston, MA, 2008.

\bibitem{I}   M.-N. Ishida, {\it Torus embeddings and de Rham complexes}, in Commutative algebra and combinatorics (Kyoto, 1985),  111--145, Adv. Stud. Pure Math., 11, North-Holland, Amsterdam, 1987. 

\bibitem{KS} H. Kim, S. Venkatesh, {\it The intersection cohomology Hodge module of toric varieties}, arXiv:2404.04767.

\bibitem{LH} J. Lipmann, M. Hashimoto, {\it Foundations of Grothendieck Duality for Diagrams of Schemes}, Lecture
Notes in Mathematics, Vol. 1960. Springer-Verlag, 2009.

\bibitem{MS} L. Maxim, J. Sch\"urmann, {\it Characteristic classes of singular toric varieties}, Comm. Pure Appl. Math. 68
(2015), no. 12, 2177--2236.

\bibitem{MS2}  L. Maxim, J. Sch\"urmann, {\it Weighted Ehrhart theory via mixed Hodge modules on toric varieties}, 
 Int. Math. Res. Not. IMRN 2025, no. 7, Paper No. rnaf067, 25 pp.

\bibitem{MSS} L. Maxim, M. Saito, J. Sch\"urmann, {\it Hirzebruch--Milnor classes of complete intersections}, 
Adv. Math. 241 (2013), 220--245.

\bibitem{Pro}  N. Proudfoot, {\it The algebraic geometry of Kazhdan–Lusztig–Stanley polynomials},
EMS Surv. Math. Sci. 5 (2018), 99--127.

\bibitem{Sa1} M. Saito, {\it Modules de Hodge polarisables}, Publ. Res. Inst. Math. Sci. 24 (1988), no. 6, 849--995.

\bibitem{Sa2} M. Saito, {\it Mixed Hodge modules}, Publ. Res. Inst. Math. Sci. 26 (1990), 221--333.

\bibitem{Sa3} M. Saito, {\it Mixed Hodge complexes on algebraic varieties}, Math. Ann. 316 (2000), 283--331.

\bibitem{Sa} M. Saito,  {\it Intersection complexes of toric varieties and mixed Hodge modules}, arXiv:2006.04081.

\bibitem{Sc} J. Sch\"urmann, {\it Characteristic classes of mixed Hodge modules}, in Topology of stratified spaces, 419--470, Math. Sci. Res. Inst. Publ., 58, Cambridge Univ. Press, Cambridge, 2011.

\bibitem{St} R. Stanley, 
{\it Generalized H-vectors, intersection cohomology of toric varieties, and related results}, in Commutative algebra and combinatorics (Kyoto, 1985), 187--213, 
Adv. Stud. Pure Math., 11, North-Holland, Amsterdam, 1987.

\bibitem{St2} R. Stanley, {\it Enumerative Combinatorics}, Volume1, Cambridge Studies in Advanced Mathematics, vol. 49. Cambridge University Press, Cambridge (2012).

\bibitem{Su} H. Sumihiro, {\it Equivariant completion}, 
J. Math. Kyoto Univ. 14 (1974), 1--28. 

\bibitem{Ta} T. Tanisaki, {\it Hodge modules, equivariant K-theory and Hecke algebras}, Publ. Res. Inst. Math. Sci. 23 (1987), no. 5, 841--879.

\bibitem{We} A. Weber, {\it Equivariant Hirzebruch class for singular varieties}, Selecta Math. (N.S.) 22 (2016), no. 3, 1413--1454.

\end{thebibliography}
\end{document}